\documentclass{amsart}

\usepackage{amsthm,amsfonts,amsmath,amssymb,latexsym,epsfig}
\usepackage{upref,amssymb,eucal}

\begin{document}
\title{The ring of differential Fourier expansions}
\author{Alexandru Buium and Arnab Saha}
\def \ep{\heartsuit}
\def \tSk{S_{k \ep}}
\def \tSp{S_{p \ep}}
\def \tSpk{S_{p_k \ep}}
\def \tSP{S_{\cP\ep}}
\def \tMk{M_{k\ep}}
\def \tMp{M_{p\ep}}
\def \tMpk{M_{p_k\ep}}
\def \tMP{M_{\cP\ep}}
\def \tIP{I_{\cP\ep}}
\def \tMunu{M_{1\ep}}
\def \Sigmat{Q}
\def \stuff{\widetilde{\tS^{\infty}}}
\def \stuffk{\widetilde{\tSk^{\infty}}}
\def \bW{{\mathbb W}}
\def \bD{{\mathbb D}}
\def \ZN{\bZ[1/N,\zeta_N]}
\def \tcS{S_{\ep}}
\def \Z{{\mathbb Z}}
\def \cI{I}
\def \cU{\mathcal U}
\def \cF{\mathcal F}
\def \tcM{M_{\Sigma}}
\def \cI{I}
\def \tcI{I_{\ep}}
\def \cK{\mathcal K}
\def \cS{S}
\def \cD{\mathcal D}
\def \cE{E}
\def \cP{\mathcal P}
\def \cA{A}
\def \cV{\mathcal V}
\def \cM{M}
\def \cN{\mathcal N}
\def \cG{\mathcal G}
\def \cB{\mathcal B}
\def \cJ{\mathcal J}
\def \tG{\tilde{G}}
\def \cF{\mathcal F}
\def \h{\hat{\ }}
\def \hp{\hat{\ }}
\def \tS{\tilde{S}}
\def \tP{\tilde{P}}
\def \tA{\tilde{A}}
\def \tX{\tilde{X}}
\def \tT{\tilde{T}}
\def \tE{\tilde{E}}
\def \tV{\tilde{V}}
\def \tC{\tilde{C}}
\def \tI{I_{\dagger}}
\def \tU{\tilde{U}}
\def \tG{\tilde{G}}
\def \tu{\tilde{u}}
\def \tx{\tilde{x}}
\def \tL{\tilde{L}}
\def \tY{\tilde{Y}}
\def \d{\delta}
\def \bZ{{\mathbb Z}}
\def \bV{{\mathbb V}}
\def \bF{{\bf F}}
\def \bE{{\bf E}}
\def \bC{{\bf C}}
\def \bO{{\bf O}}
\def \bR{{\bf R}}
\def \bA{{\mathbb A}}
\def \bZg{\bZ_{\geq 0}}
\def \bB{{\bf B}}
\def \cO{\mathcal O}
\def \ra{\rightarrow}
\def \bX{{\bf X}}
\def \bH{{\bf H}}
\def \bS{{\bf S}}
\def \bF{{\mathbb F}}
\def \bN{{\bf N}}
\def \bK{{\bf K}}
\def \bE{{\bf E}}
\def \bB{{\bf B}}
\def \bQ{{\bf Q}}
\def \bd{{\bf d}}
\def \bY{{\bf Y}}
\def \bU{{\bf U}}
\def \bL{{\bf L}}
\def \bQ{{\mathbb Q}}
\def \bP{{\bf P}}
\def \bR{{\bf R}}
\def \bC{{\mathbb C}}
\def \bM{{\mathbb M}}
\def \bG{{\mathbb G}}
\def \bP{{\bf P}}
\def \tM{M_{\ep}}
\def \tS{S_{\ep}}

\newtheorem{THM}{{\!}}[section]
\newtheorem{THMX}{{\!}}
\renewcommand{\theTHMX}{}
\newtheorem{theorem}{Theorem}[section]
\newtheorem{corollary}[theorem]{Corollary}
\newtheorem{lemma}[theorem]{Lemma}
\newtheorem{proposition}[theorem]{Proposition}
\theoremstyle{definition}
\newtheorem{definition}[theorem]{Definition}
\theoremstyle{remark}
\newtheorem{remark}[theorem]{Remark}
\theoremstyle{assumption}
\newtheorem{assumption}[theorem]{\bf Assumption}
\newtheorem{example}[theorem]{\bf Example}
\numberwithin{equation}{section}
\address{Department of Mathematics and Statistics \\ University of New Mexico \\ Albuquerque, NM 87131, USA}
\email{buium@math.unm.edu, arnab@math.unm.edu} \subjclass[2000]{11 F 32, 11 F 85}
\maketitle

\begin{abstract}
For a fixed prime we prove
structure theorems for the kernel and the image of the map that attaches to any
differential modular function its differential Fourier expansion. The image of this map, which is the ring of differential Fourier expansions, plays the role of ring of functions
on  a ``differential  Igusa curve".
Our constructions  are then
 used to perform an analytic continuation between isogeny covariant differential modular forms
 on the differential Igusa curves belonging to different primes.
\end{abstract}

\section{Introduction}

\subsection{Background}
\label{maimuta}
 The concept of {\it differential modular form} ({\it $\d$-modular form} for short) was
introduced in \cite{difmod} and further developed and applied in subsequent work, in particular in \cite{Barcau, shimura, book, eigen, local}. The present work is a continuation of this study;
 however, for the convenience of the reader,
we will not  assume here
 familiarity with the above cited papers. Indeed,
 for the purpose of this Introduction we shall begin with an informal discussion of the main concepts of this theory while later, in the body of the paper,  we shall provide a quick, yet formal, self-contained review of the necessary background.

We start by fixing a prime $p \geq 5$ and considering the ring $R:=\hat{\bZ}^{ur}_p$ obtained by completing the maximum unramified extension of the ring
of $p$-adic integers. Let $\phi:R \ra R$ be the unique lift of the $p$-power Frobenius
on $k:=R/pR$, and let  $\d_p=\d:R \ra R$ be the {\it Fermat quotient operator}
defined by
\begin{equation}
\label{ferm}
\d x:=\frac{\phi(x)-x^p}{p}\end{equation}
which, following \cite{char,book}, we view as a substitute for a ``derivative operator with respect to $p$".
Let $V$ be an affine smooth scheme over $R$ and fix a closed embedding $V \subset \bA^m$ into an affine space over $R$.
Then a map $f:V(R)\ra R$ is called a {\it $\d$-function
of order $r$} \cite{char} if
there exists a restricted power series $\Phi$ in $m(r+1)$ variables, with $R$-coefficients such that
$$f(a)=\Phi(a,\d a,...,\d^r a),$$
for all $a\in V(R)\subset R^m$.
 (Recall that {\it restricted} means {\it with coefficients converging $p$-adically to $0$}; also the definition above does not depend on the embedding $V \subset \bA^m$.)

Let $X_1(N)$ be the modular curve of level $\Gamma_1(N)$ over $R$ with $N$ not divisible by $p$; cf.\cite{DI}. In \cite{difmod} we considered the level one situation $N=1$ but here
(as in \cite{shimura,book,eigen}) we will assume $N>3$. Let $X$ be an affine open set of $X_1(N)$ disjoint from the cusps,
let $L$ be the line bundle on $X$, direct image of the sheaf of relative differentials on the universal elliptic curve over $X$, and let
$$V=Spec\left(\bigoplus_{n \in \bZ}L^{\otimes n}\right)\ra X$$
 be the  $\bG_m$-torsor associated to $L$.
 Then a
{\it $\d-$modular function}  of order $r$ and level $\Gamma_1(N)$ (holomorphic on $X$) is, by definition,  a $\d$-function
$f:V(R)\ra R$ of order $r$.

 Let $W:=\bZ[\phi]$ be the ring generated by $\phi$. For $w=\sum a_i\phi^i\in W$ ($a_i \in \bZ$)
 set $deg(w)=\sum a_i\in \bZ$; for $\lambda\in R^{\times}$ we set
 $\lambda^w:=\prod \phi^i(\lambda)^{a_i}$.
 A {\it $\d$-modular form of weight $w$} is a $\d$-modular function $f:V(R)\ra R$ such that
 $$f(\lambda \cdot a)=\lambda^w f(a),$$
  for all $\lambda\in R^{\times}$ and $a\in V(R)$, where $(\lambda,a)\mapsto \lambda \cdot a$ is the natural action $R^{\times}\times V(R)\ra V(R)$.

 We will  assume in this Introduction that the reduction mod $p$ of $X$ is contained in the ordinary locus of the modular curve.
 We  denote by $M^{\infty}$ the ring of all $\d$-modular functions and by $S^{\infty}$ the ring of modular forms of weight $0$.
There exists a natural {\it $\d$-Fourier expansion map} $$M^{\infty}\ra S^{\infty}_{for},$$
where $S^{\infty}_{for}$ is the ring of {\it $\d$-power series}, direct limit of the rings
 $R((q))[q',...,q^{(n)}]\hat{\ }$,
 where $q,q',...$ are variables and
$\hat{\ }$ means $p$-adic completion.  We may also consider the composition
 $$M^{\infty} \ra S^{\infty}_{for}\stackrel{\pi}{\ra} \widehat{S_{for}},$$
 where $S_{for}:=R((q))$ and the map $\pi$ sends $q',q'',...$ into $0$;
 we refer to this composition as the {\it Fourier expansion map}.

 One of the main  features  of this theory  \cite{difmod, shimura, book}  is   that there exist $\d$-modular forms which
possess a remarkable covariance property with respect to isogenies. These forms were called in \cite{difmod, shimura, book} {\it isogeny covariant $\d$-modular forms} and have no analogue in the classical \cite{DI} or $p$-adic \cite{Katz} theory of modular forms.
 The ring spanned by the isogeny covariant forms  is generated  by two fundamental forms $f^{\partial}$ and $f^1$ of weight $\phi-1$ and $-\phi-1$ respectively
\cite{Barcau, book}; this ring   can be viewed, in a sense explained in \cite{Barcau, book}, as
the ``projective coordinate ring" for the ``quotient of the modular
curve by the Hecke correspondences". This quotient does not exist,
of course, in usual algebraic geometry but, rather, in $\d$-{\it
geometry}; cf. \cite{book}.

As shown in \cite{difmod,shimura,book}, a fundamental role is played by the isogeny covariant $\d$-modular forms of weights $w$ with $deg(w)=-2$.
If $\d$-modular forms are morally viewed as
{\it not necessarily linear}  ``arithmetic differential operators"
(on certain line bundles over modular curves) then isogeny covariant $\d$-modular forms of weight $w$ of degree $deg(w)=-2$ should correspond to the {\it linear} ``arithmetic differential operators". So, morally,  for $\d$-modular forms,
$$(\text{isogeny covariance})+(\text{weight of degree \ $-2$)}\ \ \Longleftrightarrow \ \ \  (\text{linearity}).$$
Hence $f^1$ is ``linear" whereas $f^{\partial}$ is not.

\subsection{Aim of the paper}
\label{mohave}
The aim of the paper is two-fold namely:

1) For a fixed prime, we prove a series of  results about the  kernel and the image of the
$\d$-Fourier (respectively Fourier) expansion map; these settle some central issues left open in \cite{difmod}
and   will permit us, in particular, to introduce and study
the ``ring of functions on the $\d$-Igusa curve" and the ring of ``Igusa $\d$-modular functions".

2) We develop a ``partial differential" analogue of the above theory relative to a set of primes $\cP=\{p_1,...,p_d\}$; this
 will allow us to introduce and determine all the  ``linear partial differential operators" on (the appropriate bundles over) ``$\d_{\cP}$-Igusa curves". (No such operators exist on the appropriate bundles over  the modular curves themselves; this was indeed our main motivation for introducing $\d$-Igusa curves.)

\subsection{The theory for one prime}
\label{arizzona}
Before explaining our main results let us recall from \cite{gross}, p. 462,
the classical picture of modular forms mod $p$ (of which the ``$\d$-picture" is an analogue).
 Let $\overline{X}$ be the ordinary locus of the modular curve $X_1(N)\otimes k$  over the field $k$($=$algebraic closure of $\bF_p$), let $\overline{S}$ be
 the affine ring of $\overline{X}$, and let $\overline{M}$ be the ring of {\it modular functions} on $\overline{X}$ (by which we mean here the $k$-algebra generated by all modular forms over $k$ together with the inverse of the Hasse invariant $\overline{H}$).
  In particular, $\overline{S}\subset \overline{M}$. Furthermore let $\overline{S_{for}}=k((q))$ be the ring of Laurent power series over $k$.
 Consider the Fourier expansion map $\overline{M}\ra \overline{S_{for}}$.
This map is not injective (although its restriction to $\overline{S}$ is injective); by a Theorem
of Swinnerton-Dyer and Serre the kernel of this map is generated by $\overline{H}-1$.
 Let $\overline{\tS}$ be the image of $\overline{M}\ra \overline{S_{for}}$. This is the {\it ring of Fourier expansions} and it turns out to be a $(\bZ/p\bZ)^{\times}$-extension of $\overline{S}$; the spectrum of  $\overline{\tS}$  is birationally equivalent to the
Igusa curve. Moreover the ring $\overline{M_{\ep}}:=\overline{\tS}\otimes_{\overline{S}} \overline{M}$ corresponds, birationally, to the
 appropriate ring of modular functions on the Igusa curve.

Our   main idea in the first part of the paper
is to imitate  the above construction with
$$\overline{S},\ \ \overline{M},\ \ \overline{S_{for}}$$
 replaced by
the  rings
$$S^{\infty},\ \ M^{\infty},\ \ S^{\infty}_{for},$$
 where  $M^{\infty}$ is the ring  of $\d$-modular functions,  $S^{\infty}$ is the  ring of $\d$-modular forms of weight $0$, and  $S^{\infty}_{for}$ is the ring of {\it $\d$-power series}; cf. section \ref{maimuta}.
  Let $\tS^{\infty}$ be the image of the $\d$-Fourier expansion map $M^{\infty}\ra S^{\infty}_{for}$. This is the {\it ring of $\d$-Fourier expansions} and
will be viewed as playing  the role of ring of functions on a ``$\d$-Igusa curve". Moreover a certain ``partially completed version", $M^{\infty}_{\ep}$, of the ring
$\tS^{\infty}\otimes_{S^{\infty}} M^{\infty}$  will play
 the role  of  ring of ``Igusa $\d$-modular functions". (We will not introduce, in this paper, an object called the {\it $\d$-Igusa curve}; such an object can be formally introduced in the sense of $\d$-algebraic geometry \cite{book} but we shall not pursue this here. All we shall be working with are
certain rings that play the roles of rings of functions, or rings of  sections of bundles, on such a $\d$-geometric object.)

Here are (somewhat rough formulations of) our main results about the kernel and the image of the $\d$-Fourier (respectively Fourier) expansion map. For  some of the  terminology involved,
and for more precise formulations of the results we refer to the body of the paper, as we shall explain presently.
The first two Theorems below
should be viewed as  $\d$-analogues  of the Swinnerton-Dyer and Serre Theorem about the kernel
 of the Fourier expansion map in positive characteristic.

\begin{theorem}
\label{teo2} The kernel of the $\d$-Fourier expansion map, $M^{\infty}\ra S^{\infty}_{for}$,
is the $p$-adic closure of the ideal generated by the elements $\d^i(f^{\partial}-1)$, where $i \geq 0$.
\end{theorem}

\begin{theorem}
\label{teo3}  The kernel of the Fourier expansion map,
$M^{\infty} \ra \widehat{S_{for}}$,
 is the $p$-adic closure of the ideal  generated by the elements $\d^i(f^{\partial}-1)$ and $\d^i f^1$, where $i\geq 0$.\end{theorem}

In the next statement, for any ring $A$, we denote by $\overline{A}$ the ring $A/pA$.

\begin{theorem}
\label{teo1}
The cokernel of the $\d$-Fourier expansion map  $M^{\infty} \ra S^{\infty}_{for}$ is torsion free.
In particular the ring $\overline{\tS^{\infty}}$ is an integral domain. Moreover $\overline{\tS^{\infty}}$ is an ind-\'{e}tale $\bZ_p^{\times}$-extension  of
 $\overline{S^{\infty}}$.
\end{theorem}

The statement about the cokernel  should be viewed as a  {\it $\d$-expansion principle} for Igusa $\d$-modular forms.
 The rest of the Theorem  shows that, morally,
the  the ``$\d$-Igusa curve" is irreducible and is
a  pro-\'{e}tale ``formal $\bZ_p^{\times}$-cover" of the ``$\d$-modular curve'' (i.e. of the modular curve viewed as an object of ``$\d$-geometry'').
 By the way, we will also show,  in the body of the paper, that the ring $\tS^{\infty}$   comes
 equipped with a  sequence of derivations naturally associated to  the Serre derivation operator; these derivations ``topologically generate"
  the space of all
derivations on  $\tS^{\infty}$.
Theorems \ref{teo1}, \ref{teo2}, \ref{teo3} are consequences of the (more precisely formulated) Theorem \ref{babanovac} in the body of the paper.

In the the following Theorem \ref{teo4}
 we assume that the reduction mod $p$ of $X$ is the whole of the ordinary locus of the modular curve.
To state this theorem recall Katz' rings
$$\bD \subset \bW \subset \widehat{S_{for}}=R((q))\h$$
  where $\bD$ is the ring of {\it divided congruences} and $\bW$ is the ring of {\it generalized $p$-adic modular forms}, both with coefficients in $R$, and both viewed as embedded into $R((q))\h$ via the Fourier expansion. Cf. \cite{gouvea, Katzgen}, and also the review in the present paper. Recall that $\bD\subset R[[q]]$; also, if $\Delta \in R((q))$ is the discriminant  then
$\bD+R[\Delta^{-1}]$ is $p$-adically dense in $\bW$.

 \begin{theorem}
\label{teo4}
 The  image of the Fourier expansion map
 $M^{\infty} \ra \widehat{S_{for}}$ contains $\bD$ and hence is  $p$-adically dense in   $\bW$.
 \end{theorem}

Morally  Theorem \ref{teo4} (which is Corollary \ref{umbirac} in the body of the paper)
exhibits the ``world of $\d$-modular functions" as a lift (with ``huge" kernel described in  Theorem \ref{teo3})
 of the ``world of generalized $p$-adic modular functions" of Katz.
 One can then ask if some of the basic constructions with values in Katz' ring $\bW$ (such as  the $U$-operator, various measures, Galois representations, etc.) can be lifted
 naturally to the world of $\d$-modular functions.

  Note that Theorem \ref{teo4} has the following  consequence that
  is independent of our theory.
 Let $\d_0:R((q))\hat{\ }\ra R((q))\hat{\ }$ be the operator
$$\d_0(\sum a_n q^n):=\frac{\sum \phi(a_n)q^{np} -(\sum a_nq^n)^p}{p},$$
 let $K=R[1/p]$,  let $M(R,\kappa,N)$ denote the space of modular forms
over $R$ of weight $\kappa$ and level $\Gamma_1(N)$, let $\Delta \in M(R,12,N)$ be the discriminant form, and let $E_{p-1}\in M(R,p-1,N)$ be the normalized  Eisenstein form of weight $p-1$. Then Theorem \ref{teo4} implies that
any series  $f(q)\in \bD$
can be represented in $R((q))\hat{\ }$ as
\begin{equation}
\label{doicur}
f(q)=\Phi(f_1(q),...,f_n(q),...,\d_0^r(f_1(q)),...,\d_0^r(f_n(q))),
\end{equation}
where $f_j=F_j\Delta^{\mu_j} E_{p-1}^{\nu_j}$,  $F_j \in M(R,\kappa_j,N)$, $\kappa_j,\mu_j,\nu_j \in \bZ$,
$f_j(q)\in R((q))$ are the Fourier expansions of $f_j$, and $\Phi$ is a restricted power series in $n(r+1)$ variables with $R$-coefficients.
So, morally,  the elements
   of Katz's ring $\bD$ of divided congruences  can be realized as limits of very special
 divided congruences that arise by iterating the Fermat quotient operation. It is not clear to us if this statement can be proved directly, independently of our theory.

\medskip

\subsection{The theory for several primes}
\label{orangutan}
As an application of the above one prime constructions
we will  introduce and study  {\it linear partial differential operators} in the setting of modular curves, with respect to  $d$ arithmetic directions; these directions are represented by a set of primes $\cP=\{p_1,...,p_d\}$ (along which the ``derivatives"
will be the corresponding Fermat quotient operators). Such a theory for algebraic groups (rather than modular curves) was developed in \cite{laplace} where arithmetic analogues of Laplacians were constructed on the additive group, the multiplicative group, and on elliptic curves over $\bQ$. The basic idea in \cite{laplace} was, very roughly speaking,  to construct {\it linear arithmetic partial differential operators} along the ``vertical divisors" corresponding to each of the
primes $p_1,...,p_d$ and then to perform a sort of {\it analytic continuation} between the various primes
 along a ``horizontal divisor".
We will keep this point of view of analytic continuation in the present paper. (For the convenience of the reader, we will make the present paper essentially independent of \cite{laplace}; but for an informal explanation of  analytic continuation, the reader may want to consult  the Introduction to \cite{laplace}). With this point of view we will be able to achieve our program of constructing linear partial differential operators in the modular setting; but for this we will have to pay  the price of passing from modular curves to   {\it $\d_{\cP}$-Igusa curves} (a several primes generalization of the $\d$-Igusa curves).  The main reason why passing
 from modular curves to $\d_{\cP}$-Igusa curves is crucial is that, unlike the former,  the latter carry a certain ``tautological"
weight one, order zero  form, which we shall call $f^0$, and which will be the key  to performing ``analytic continuation" between various primes;
we will show that such an analytic continuation cannot be performed in the context of  the modular curves themselves.

Our main result here is, morally, a complete determination of all {\it linear arithmetic partial differential operators}  on the
appropriate bundles over the ``$\d_{\cP}$-Igusa curves". The technical way to express this is the following (roughly formulated):

\begin{theorem}
Let $w$ be a weight of degree $deg(w)=-2$.
Then the module of all  {\it weight $w$ isogeny covariant Igusa $\d_{\cP}$-modular forms} of  weight $w$ and order $r=(r_1,...,r_d)$
is free of rank $r_1r_2...r_d$.
\end{theorem}

Cf. Theorem \ref{uniqueness} in the body of the paper (and the discussion preceding it) for a precise formulation. We mention that, intuitively,
``order $r=(r_1,...,r_d)$" means ``order $r_k$ with respect to $p_k$ for each $k=1,...,d$".

The forms in the above Theorem
are constructed using the tautological form $f^0$ and another basic form $f^e$ of weight $0$
and order $e:=(1,...,1)$. The form $f^e$ itself is constructed using the forms
 $f^1$ corresponding to the  various primes and can be viewed
 as  a modular analogue of the arithmetic Laplacians in \cite{laplace}.

By the way we will also determine all isogeny covariant Igusa $\d_{\cP}$-modular forms of any order and any weight $w$ with $deg(w)=0$; they are all obtained from $f^0$. Finally we shall be able to
 analytically continue the differential eigenforms of \cite{eigen}; the resulting form will be called $f^{2e}$, will have weight $0$ and order $2e=(2,...,2)$, will {\it not} be isogeny covariant, and will  live on the modular curves themselves (rather than on the ``$\d_{\cP}$-Igusa curves").

\subsection{Variants}

In all the discussion above, we considered the  {\it modular curves}  $X_1(N)$ parameterizing elliptic curves with $\Gamma_1(N)$ level structures; nevertheless a substantial part of the  theory will be developed also  in the case of  Shimura curves $X^D(\cU)$ parameterizing false elliptic curves with {\it level $\cU$   structures} \cite{buzzard}. The role of Fourier expansions in the Shimura curve case will be played by Serre-Tate expansions,
in the sense of \cite{mori,shimura,book}. We will also consider the theory
over modular curves with respect to Serre-Tate expansions.

Finally note that our theory in $d$ arithmetic dimensions can be viewed as a $0+d$ dimensional theory (where $0$ indicates that we have no ``geometric direction'' on the base).
On the other hand we developed in \cite{pdemod} a $1+1$ dimensional theory (where the base has one geometric and one arithmetic direction.)
One  can then ask if the $0+d$ dimensional theory in the present paper and the $1+1$ dimensional theory in \cite{pdemod}
can be ``unified'' as parts of a $1+d$ dimensional theory. At this point it is not clear whether this is possible; cf. the last section of \cite{laplace}
for comments on the difficulties arising from such an attempt at  unification.

\subsection{Analytic analogues}
One can ask if the basic forms $f^e$, $f^{2e}$ referred to in section \ref{orangutan}
 have analogues in (real/complex) analysis.
The forms $f^{2e}$ are intimately related
to the arithmetic Laplacians in \cite{laplace} whose analytic analogues are discussed in the Introduction of that paper.
The form $f^e$, on the other hand, can be loosely viewed as having an analytic analogue which we now describe.

Let $D$ be a domain in the complex $z$-plane. Let ${\mathbb H}:=\{\tau \in \bC;Im(\tau)>0\}$
be the upper half plane. Let $F \subset {\mathbb H}$ be a domain containing none of the fixed points
of the $SL_2(\bZ)$-action, such that any two $SL_2(\bZ)$-conjugate points in $F$ are conjugate under a translation by an integer.
(E.g. one can take $F=\{\tau \in \bC;Im(\tau)>1\}$.)
Let $j:{\mathbb H} \ra \bC$ be the classical $j$-function, let $G=j(F)$, and let $j^{-1}:G \ra F$ be the (multivalued) inverse of $j$.
For any $C^{\infty}$-function $u \in C^{\infty}(D,G)$ we denote by $q_u \in C^{\infty}(D,\bC^{\times})$ the (well defined !) function
$$q_u(z)=e^{2 \pi \sqrt{-1} \cdot j^{-1}(u(z))}.$$
Then our form $f^e$ can be viewed as an arithmetic analogue of the Laplace-type operator
$$C^{\infty}(D,G) \ra C^{\infty}(D,\bC)$$ defined by
$$u \mapsto \partial_z \partial_{\bar{z}} \log q_u=\partial_z \left( \frac{\partial_{\bar{z}} q_u}{q_u}\right)=
\partial_{\bar{z}} \left( \frac{\partial_{z} q_u}{q_u}\right).$$
The last equalities can be viewed as a {\it Dirac decomposition} for our operator and this decomposition will have an arithmetic analogue in the case of $f^e$.
Cf. \cite{laplace} for more on Dirac decompositions.

 \subsection{Plan of the paper}
 The paper has two parts. In the first part we develop the theory for one prime, first in
 an axiomatic setting, for
  an arbitrary  curve equipped with an arbitrary  line bundle, and then in the concrete setting of modular or Shimura curves equipped with their bundles of modular forms.
  The main result here is Theorem \ref{babanovac}  (which implies Theorems \ref{teo1}, \ref{teo2}, \ref{teo3}) and Corollary
\ref{umbirac} (whose content is that of Theorem \ref{teo4}).

   Our one prime constructions will be used in the second part of the paper, where the theory for at least two primes is developed. In this second part analytic continuation in the modular/Shimura context is introduced and the main results on the space of isogeny covariant Igusa $\d$-modular forms,
   referred to in section \ref{orangutan}, are stated and proved. Cf. Theorems \ref{uniqueness}, \ref{buniqueness}.

\subsection{Acknowledgment} While writing this paper the first author was
partially supported by NSF grants DMS-0552314 and DMS-0852591.  Any
opinions, findings, and conclusions or recommendations expressed in
this material are those of the author and do not necessarily reflect
the views of the National Science Foundation.

\section{The theory for one prime}
\subsection{Review of concepts and terminology from \cite{char, book}}
Unless otherwise stated all rings and algebras will be commutative with unit element.
For any $A$-algebra $\varphi:A \ra B$ and any element $a \in A$ we continue to denote by $a$ the element $\varphi(a)=a 1_B$.
We fix, throughout this paper, a prime integer
$p \geq 5$. For any $\bZ$-module $M$ we set $M\h=\varprojlim M/p^nM$,
the $p$-adic completion of $M$ and $\overline{M}:=M/pM=M \otimes \bZ/p\bZ$, the reduction of $M$ mod $p$.
We say $M$ is $p$-adically complete if $M \ra M\hat{\ }$ is an isomorphism.
For $m \in M$ we let $\overline{m} \in \overline{M}$ be the image of $m$.
For any scheme $X$ we set $\overline{X}:=X \otimes \bZ/p\bZ$. We denote by $\bZ_{(p)}$ the local ring of $\bZ$ at $(p)$.
We will repeatedly use the fact that if $M \ra N$ is a homomorphism of $\bZ_{(p)}$-modules
such that $\overline{M}\ra \overline{N}$ is injective, $p$ is a non-zero divisor in $N$, and $M$ is $p$-adically separated then $M \ra N$ is injective and has torsion free cokernel; and that,
conversely, if $M \ra N$ is a morphism of $\bZ_{(p)}$-modules which is injective and has torsion free cokernel then $\overline{M}\ra \overline{N}$ is injective.

\subsubsection{$p$-derivations}
Let $C_p(X,Y) \in \bZ[X,Y]$ be the polynomial with
integer coefficients \[C_p(X,Y):=\frac{X^p+Y^p-(X+Y)^p}{p}.\] A
$p-${\it derivation} from a ring $A$ into an $A-$algebra
$\varphi:A \ra B$ is a map $\d:A \ra B$ such that $\d(1)=0$ and
\[\begin{array}{rcl}
\d(x+y) & = &  \d x + \d y
+C_p(x,y)\\
\d(xy) & = & x^p \cdot \d y +y^p \cdot \d x
+p \cdot \d x \cdot \d y,
\end{array}\] for all $x,y \in A$. Given a
$p-$derivation we always denote by $\phi:A \ra B$ the map
$\phi(x)=\varphi(x)^p+p \d x$; then $\phi$ is a ring homomorphism. A
{\it prolongation sequence} is a sequence $S^*=(S^n)_{n \geq 0}$ of  rings $S^n$, $n
\geq 0$, together with ring homomorphisms $\varphi_n:S^n \ra
S^{n+1}$ and $p-$derivations $\d_n:S^n \ra S^{n+1}$ such that
$\d_{n+1} \circ \varphi_n=\varphi_{n+1} \circ \d_n$ for all $n$. We
usually denote all $\varphi_n$ by $\varphi$ and all $\d_n$ by $\d$
and we view $S^{n+1}$ as an $S^n-$algebra via $\varphi$. A morphism
of prolongation sequences, $u^*:S^* \ra \tilde{S}^*$ is a sequence
$u^n:S^n \ra \tilde{S}^n$ of ring homomorphisms such that $\delta
\circ u^n=u^{n+1} \circ \d$ and $\varphi \circ u^n=u^{n+1} \circ
\varphi$. Let $W$ be the ring of polynomials $\bZ[\phi]$ in the
indeterminate $\phi$. Then, for $w=\sum_{i=0}^r a_i \phi^i \in W$,
we set $deg(w):=\sum a_i$. If $a_r \neq 0$ we set $ord(w)=r$; we
also set $ord(0)=0$. For $w$ as above (respectively for $w \in
W_+:=\{\sum b_i \phi^i\ |\ \ b_i \geq 0\}$), $S^*$ a prolongation
sequence, and $x \in (S^0)^{\times}$ (respectively $x \in S^0$) we
can consider the element $x^w:=\prod_{i=0}^r \varphi^{r-i}
\phi^i(x)^{a_i} \in (S^r)^{\times}$ (respectively $x^w \in S^r$). We
let $W(r):=\{w \in W\ |\ ord(w) \leq r\}$.

Let $R:=R_p:=\hat{\bZ}_p^{ur}$ be the completion of the maximum
unramified extension of the ring of $p$-adic integers $\bZ_p=(\bZ_{(p)})\h$
 and we denote by $k$ its residue field, $k=R/pR$. Then $R$ has a unique
$p-$derivation $\d:R \ra R$ given by $$\d x=(\phi(x)-x^p)/p,$$ where
$\phi:R \ra R$ is the unique lift of the $p-$power Frobenius map
on $k$. One can consider the  prolongation sequence $R^*$ where
$R^n=R$ for all $n$. By a {\it prolongation sequence over $R$} we
understand a prolongation sequence $S^*$ equipped with a morphism
$R^* \ra S^*$. From now on all our prolongation sequences are
assumed to be over $R$.

By a $\d$-ring we mean a ring together with a $p$-derivation on it. A morphism of $\d$-rings
is a ring homomorphism that commutes with the given $p$-derivations. In what follows all $\d$-rings
will be assumed {\it over $R$} (i.e. equipped with $\d$-ring homomorphisms from $R$.)
If $a$ is an element of a $\d$-ring we will sometimes denote by $a',a'',...,a^{(r)}$ the sequence $\d a, \d^2 a,...,\d^r a$.

\subsubsection{Conjugate derivations}
Let $A$ be a $\d$-ring in which $p$ is a non-zero divisor and let $u:A^0\ra A$ be a ring
homomorphism. Let $\partial:A^0\ra A^0$ be an $R$-derivation and let $j \geq 0$ be an integer.
An $R$-derivation $\partial_j:A \ra A$ will be called a {\it $j$-conjugate} of $\partial$
if for any integer $s \geq 0$ we have
$$\partial_j \circ \phi^s \circ u=\delta_{js} \cdot p^j \cdot \phi^s \circ u \circ \partial:A^0\ra A,$$
where $\delta_{js}$ is the Kronecker symbol. A sequence $(\partial_j)_{j \geq 0}$ where for  each $j$, $\partial_j:A \ra A$
is a  $j$-conjugate of $\partial$, will be referred to as a {\it complete sequence of conjugates} of $\partial$.
Let us say that $A$ is {\it topologically $\d$-generated} by $A^0$ if the smallest $\d$-subring of $A$ that contains $u(A^0)$
is $p$-adically dense in $A$. It is then trivial to see that if $A$ is $p$-adically separated and topologically $\d$-generated by $A^0$
then 1) any derivation $\partial:A^0\ra A^0$ has at most one $j$-conjugate $\partial_j:A\ra A$ for each $j$ and 2)
if $(\partial_j)_{j \geq 0}$ is a complete sequence of conjugates of $\partial$ then, for all $j \geq 0$,
$$\begin{array}{rcl}
\partial_j \circ \phi & = & p\cdot \phi \circ \partial_{j-1}:A \ra A\\
\  & \  & \  \\
\partial_j \circ \d^s \circ u & = & 0:A^0\ra A,\ \text{for}\ \ s<j,\\
\  & \  & \  \\
\partial_j \circ \d^j \circ u & = & \phi^j\circ u \circ \partial:A^0\ra A.\end{array}$$
Here $\partial_{-1}=0$.

\subsubsection{$p$-jet spaces}
Given a scheme $X$ of finite type over $R$ we introduced in \cite{char} a sequence of formal ($p$-adic) schemes over $R$, called the {\it $p$-jet spaces of $X$}, which we denoted by $J^r(X)$, $r \geq 0$. In case $X$ is affine,
$X=Spec\ R[x]/(f)$, with $x$ a tuple of indeterminates and $f$ a tuple of polynomials, we have
$$J^r(X)=Spf\ R[x,x',...,x^{(r)}]\h/(f, \d f,...,\d^r f)$$
where $x',...,x^{(r)}$ are new tuples of variables and $R[x,x',...,x^{(r)}]\h$ is a prolongation sequence via $\d x=x'$, $\d x'=x''$,... For $X$ not necessarily affine we set $\cO^r(X):=\cO(J^r(X))$; these rings form a prolongation sequence.
 If $X$ is affine the prolongation sequence $(\cO^r(X))_{r \geq 0}$ has the following universality property: if $(S^r)_{r \geq 0}$ is any
 prolongation sequence over $R$ of $p$-adically complete rings $S^r$ and $u:\cO(X)\ra S^0$ is any $R$-algebra homomorphism then there exists a unique morphism of prolongation sequences over $R$, $u^r:\cO^r(X)\ra S^r$, such that $u^0$ induces $u$. By this universality property, for $X$ not necessarily affine,
 each element of $\cO^r(X)$ naturally  defines a function $X(R)\ra R$. Such functions are called in \cite{char,book}  {\it $\d$-functions of order $r$}. If $X$ is smooth then any element of $\cO^r(X)$ is uniquely determined by the induced $\d$-function $X(R)\ra R$.
We set $\cO^{\infty}(X):=\lim_{\ra} \cO^r(X)$. If $X/R$ is smooth with $\overline{X}$ connected
then the schemes $\overline{J^r(X)}$ are smooth varieties over $k$. Moreover  $\cO^r(X)$, $\cO^{\infty}(X)$ are integral domains, and
 $p$ is a prime element in these rings.
In addition, $\cO^{\infty}(X)$ is $p$-adically separated and topologically $\d$-generated by $\cO(X)$ (and hence also by $\cO^0(X)=\cO(X)\h$).
If $X \ra Y$ is an \'{e}tale morphism then $J^r(X)\simeq J^r(Y)\widehat{\times}_{\widehat{Y}} \widehat{X}$.

Recall from \cite{book}, Proposition 3.45, that if $X/R$ is smooth  then for any $R$-derivation $\partial:\cO^0(X)\ra \cO^0(X)$ there exists
a (necessarily unique) complete sequence of conjugates $\partial_j:\cO^{\infty}(X) \ra \cO^{\infty}(X)$ of $\partial$. Moreover $\partial_j \cO^s(X) \subset \cO^s(X)$ for all $j,s \geq 0$ and $\partial_j \cO^s(X)=0$ for $s<j$.

\subsection{The axiomatic theory}

In this section we develop the one prime version of the theory of this paper in the axiomatic setting of an arbitrary curve equipped
with a line bundle. In the next section we will specialize our discussion to the case of modular (respectively Shimura) curves
and their natural bundles of modular forms. The main result of this axiomatic section is Theorem \ref{main}.
This Theorem will later be strengthened,
 in the concrete setting of modular curves; cf. Theorem \ref{babanovac}. The strengthened version  will morally say that the $\d$-Igusa curve
is a connected pro-\'{e}tale $\bZ_p^{\times}$-cover of the modular curve.

\subsubsection{Framed curves}
We start with the following data:
\begin{equation}
\label{data1}
X,\ L,
\end{equation}
where $X$ is a  smooth  affine curve over $R$ with connected reduction mod $p$,
$\overline{X}$,
and $L$ is an invertible sheaf  on $X$ which we identify with its module of global sections.
Consider the scheme
$$V:=Spec\left( \bigoplus_{n \in \bZ} L^{\otimes n}\right).$$
For any $R$-algebra $B$ the set of $B$-points $V(B)$ naturally identifies with the set of pairs $(P,\xi)$
where $P\in X(B)$ and $\xi$ is a basis of the pull back of $L$ to $Spec\ B$ by $P$. In particular $V\ra X$ is a $\bG_m$-torsor
with respect to the action $B^{\times}\times V(B)\ra V(B)$ given by
\begin{equation}
\label{joe}
(\lambda,(P,\xi))\mapsto \lambda \cdot (P,\xi):=(P,\lambda^{-1}\xi).\end{equation}
Set
\begin{equation}
\label{putoare}
\begin{array}{lll}
S & := & S_X :=\cO(X),\\
\  & \  &   \\
 M  & := & M_X  :=\cO(V)=\bigoplus_{n \in \bZ} L^{\otimes n}\end{array}\end{equation}
We assume in what follows that we are given one more piece of data, namely
either an $R$-point of $X$:
\begin{equation}
\label{data3}
P \in X(R)
\end{equation}
or an open embedding $X \subset X^*$
into a smooth curve $X^*$ over $R$, with connected reduction mod $p$, $\overline{X^*}$,  plus  an $R$-point
 of the reduced closed subscheme $X^*\backslash X$ of $X^*$,
 \begin{equation}
 \label{data5}
 P \in (X^* \backslash X)(R).
 \end{equation}
 In the first case (\ref{data3}) we set $S_{for}=R[[t]]$, where $t$ is a variable.
 In the second case (\ref{data5}) we set $S_{for}:=R((q)):=R[[q]][q^{-1}]$, where $q$ is a variable.
Assume we are given an isomorphism between $Spf\ R[[t]]$ (respectively $R[[q]]$)
and the completion of $X$
(respectively $X^*$) along the image of $P$.
Then, in both cases we have an induced homomorphism $E:S \ra S_{for}$ which is injective,
inducing an injective map $\overline{E}:\overline{S}\ra \overline{S_{for}}$.
We shall finally assume that we are given yet another piece of data, namely an extension of  $E$
to  a homomorphism
\begin{equation}
\label{data6}
\cE:M \ra S_{for}.
\end{equation}
It is convenient to collect all the above data in one definition as follows:

\begin{definition}
\label{ED}
A {\it framed curve} is a tuple $X,L,P,E$ where $X,L$ are as in (\ref{data1}), $P$ is either as in (\ref{data3}) or as in (\ref{data5}), and $E$
is as in (\ref{data6}). For simplicity we also say that $X$ is a framed curve with frame $L,P,E$. We say that $X$ is {\it Serre-Tate-framed}, respectively  {\it Fourier-framed}, according as $P$ is as in (\ref{data3}) or as in (\ref{data5}). Accordingly $E$ is referred to as a {\it Serre-Tate expansion map} respectively a {\it Fourier  expansion map}.
\end{definition}

 \subsubsection{The rings $S^r,M^r$}
 Assume we are given a framed curve $X=Spec\ S$. We may consider the following rings:
 \begin{equation}
 \begin{array}{lll}
S^r & := & S^r_X:=\cO^r(X),\ \ r\geq 0\\
\  & \  & \ \\
M^r & := & M^r_X:=\cO^r(V),\ \ \ r\geq 0\\
\ & \  & \\
S^{\infty} & := & \lim_{\ra} S^r,\\
\ & \  & \   \\
 M^{\infty} &:= &\lim_{\ra} M^r.\end{array}
\end{equation}
An element $f \in M^r$ is said to be of weight $w \in W$ if, and only if, the induced $\d$-function $f:V(R)\ra R$ satisfies
$$f(\lambda \cdot a)=\lambda^wf(a)$$
 for all $\lambda \in R^{\times}$, $a \in V(R)$, where $(\lambda,a)\mapsto \lambda \cdot a$ is the natural action $R^{\times}\times V(R)\ra V(R)$; cf. (\ref{joe}). We denote by $M^r(w)=M^r_X(w)$ the $R$-module of all elements of $M^r=M^r_X$ of weight $w$.

 If $L$ is trivial on $X$ and $x$ is a basis of $L$ then we have  identifications
$$\begin{array}{rcl}
M & = & S[x,x^{-1}],\\
\ & \  & \\
 M^r & = & S^r[x,x^{-1},x',...,x^{(r)}]\h,\\
 \  & \  & \\
M^r(w) & = & S^r \cdot x^w\subset M^r.\end{array}$$
 We may also consider the element
$x^{-1}\otimes x\in M\otimes_S M$ which we refer to as the {\it tautological element of $M\otimes_S M$}.
Clearly  $x^{-1}\otimes x$ do not depend on the choice of the basis $x$ of $L$.

By \cite{book}, Proposition 3.14, the rings  $\overline{S^r}$ are integral domains,
and the maps $\overline{S^r}\ra \overline{S^{r+1}}$ are injective. In particular the rings
$S^r$ are integral domains and the maps $S^r \ra S^{r+1}$ are injective with torsion free cokernels.
The analogous  statements hold for $M^r$.
So, in particular, $\overline{S^{\infty}}$ and $\overline{M^{\infty}}$ are integral domains.

Let $t',t'',...$ and $q',q'',...$ be new variables and consider the prolongation sequence $(S^r_{for})_{r \geq 0}$,
$$
S^r_{for}  =  R[[t]][t',...,t^{(r)}]\h,$$
respectively
$$S^r_{for}  =  R((q))\h[q',...,q^{(r)}]\h.$$
We set
$$S_{for}^{\infty}:=\lim_{\ra} S^r_{for}.$$
Then the expansion maps induce, by universality,   morphisms of prolongation sequences,
\begin{equation}
\label{adelita}
\cE^r:M^r\ra S^r_{for};\end{equation}
the maps $\cE^r$ will be referred to as $\d$-{\it expansion maps} for $M^r$.
They induce a $\d$-{\it expansion map}
\begin{equation}
\label{woody}
E^{\infty}:M^{\infty}\ra S^{\infty}_{for}.\end{equation}

We have the following {\it $\d$-expansion principle} for $S^r$:

\begin{proposition}
\label{expansionprinciple}
The induced map $$\overline{\cE^r}:\overline{S^r} \ra \overline{S^r_{for}}$$ is injective. In particular,
$\cE^r:S^r \ra S^r_{for}$, and hence  the $\d$-expansion maps $$\cE^r:M^r(w) \ra S^r_{for}$$ are injective, with torsion free cokernel.
\end{proposition}

(In this paper the words ``torsion free", without the specification ``as an $A$-module",
will always mean ``torsion free as a $\bZ$-module".)

{\it Proof.} In case (\ref{data3}) this is \cite{book}, Proposition 4.43. Case (\ref{data5}) follows easily from case (\ref{data3}) applied to $X^*$ instead of $X$.
\qed

\subsubsection{The rings $\tcS^r$}
Next, for a framed curve $X=Spec\ S$, we define
\begin{equation}
\label{tillde}
\begin{array}{rcl}
\tS^r & := & Im(\cE^r:M^r\ra S^r_{for})\\
\  & \  & \  \\
\tS^{\infty} & := & \lim_{\ra} \tS^r=Im(E^{\infty}:M^{\infty}\ra S^{\infty}_{for}). \end{array}\end{equation}
The ring $\tS^{\infty}$ will later morally play the role
of ``coordinate ring of the $\d$-Igusa curve".

The following is trivial to check (using the definitions and Proposition \ref{expansionprinciple}):

\begin{proposition}
\label{lotsofinj}
\

1) The homomorphisms
$\overline{S^r}\ra \overline{\tS^r}$, $\overline{S^{\infty}}\ra \overline{\tS^{\infty}}$
are injective.
 In particular the homomorphisms $S^r\ra \tS^r$, $S^{\infty}\ra \tS^{\infty}$ are injective with torsion free cokernel.

2) The homomorphisms $\tS^r\ra \tS^{r+1}$ are injective.
\end{proposition}

\begin{remark}
The ring $\overline{\tS^{\infty}}$ is not a priori an  integral domain and the map $\overline{\tS^{\infty}} \ra \overline{S^{\infty}_{for}}$
is not a priori injective. The
ring $\overline{\tS^{\infty}}$, however, has a natural quotient which is an integral domain, namely:
\begin{equation}
\label{doitild}
\stuff:=Im(\overline{M^{\infty}}\ra \overline{S^{\infty}_{for}}).
\end{equation}
This ring is going to play a role in what follows.
We will prove later that, in the concrete setting of modular curves the map $\overline{\tS^{\infty}}\ra \overline{S^{\infty}_{for}}$ is injective
($\d$-expansion principle)
hence
$\overline{\tS^{\infty}}$ is an integral domain and the surjection $\overline{\tS^{\infty}}\ra \stuff$
is an isomorphism. Cf. Theorem \ref{babanovac}.
\end{remark}

\begin{definition}
\label{ass2}
A framed curve is called {\it ordinary} if
there exists an element $f \in M^1(\phi-1)$ which is invertible in the ring $M^1$,
such that $\cE^1(f)=1$.
\end{definition}

Note that if an $f$ as above exists then, by Proposition \ref{expansionprinciple},  $f$ is necessarily unique.
The terminology {\it ordinary} will be justified later in our applications to modular and Shimura curves.

In what follows we will analyze in some detail the structure of the rings $\tS^{\infty}$ and the various rings constructed from it.

\begin{definition}
\label{ooyy}
Let $A$ be a $k$-algebra where $k$ is a field. Let   $A \subset B$ a ring extension, and $\Gamma$ a profinite abelian group
acting on $B$ by   $A$-automorphisms. We say that $B$ is a $\Gamma$-{\it extension} of $A$ if one can write $A$ and $B$ as filtered unions of
finitely generated $k$-subalgebras, $A=\bigcup A_i$, $B=\bigcup B_i$, indexed by some partially ordered set, with $A_i\subset B_i$,
and one can write $\Gamma$ as an inverse limit of finite abelian groups, $\Gamma=\varprojlim \Gamma_i$, such that
the $\Gamma$-action on $B$ is induced by a system of compatible $\Gamma_i$-actions on $B_i$
 and
$$B_i^{\Gamma_i}=A_i$$
for all $i$.
(Then, of course, we also have $B^{\Gamma}=A$.) If in addition one can choose the above data such that each $A_i$ is smooth over $k$ and each $B_i$ is \'{e}tale over $A_i$
we say that $B$ is an {\it ind-\'{e}tale $\Gamma$-extension} of $A$.
\end{definition}

Here are a couple of easy facts about this concept:

\begin{lemma}
\label{escr}
\

1) Assume $B$ is a $\Gamma$-extension of $A$ and $C:=B/I$ is a quotient of $B$ by an ideal $I$. Then
$C$ is integral over $A$.

2) Assume $B$ is an ind-\'{e}tale $\Gamma$-extension of $A$ and let $I$ be a prime ideal of $B$ such that $I\cap A=0$.
Then $C:=B/I$ is an ind-\'{e}tale $\Gamma'$-extension of $A$ where $\Gamma'$ is a closed subgroup of $\Gamma$.
\end{lemma}

{\it Proof}.
Assertion 1 is clear. Let's prove assertion 2. Using the notation in Definition \ref{ooyy} set $Y_i=Spec\ B_i$, $V_i:=Spec\ A_i$, $Z_i:=Spec\ C_i$, $C_i:=B_i/B_i\cap I$.
Let $\Gamma_i':=\{\gamma \in \Gamma_i;\gamma Z_i=Z_i\}$.
By Lemma \ref{fomizt} below $C_i$ is \'{e}tale over $A_i$ and $C_i^{\Gamma'_i}=A_i$ so one can take
$\Gamma':=\varprojlim \Gamma'_i$ acting on $C=\varprojlim C_i$.
\qed

We have used the following ``well known'' lemma (whose proof will be ``recalled" for convenience):

\begin{lemma}
\label{fomizt}
Let $V$ be a smooth affine variety over a field $k$, let $Y \ra V$ be a finite \'{e}tale map, and let $G$ be a finite abelian group acting
on $Y$ such that $Y/G=V$. Let $Z \subset Y$ be a subvariety that dominates $V$ and let $G'=\{\gamma \in G;\gamma Z=Z\}$.
Then $Z$ is a connected component of $Y$ (hence is \'{e}tale over $V$) and $Z/G'=V$.
\end{lemma}

{\it Proof}.
Since $V$ is smooth the connected components $Z_1,...,Z_n$ of $Y$ are irreducible so $Z$ is a connected component of $Y$, say $Z=Z_1$.
Since $V$ is connected $G$ acts transitively on the set $\{Z_1,...,Z_n\}$ hence the stabilizers in $G$ of the various $Z_i$s are conjugate in $G$,
hence they are equal, because $G$ is abelian. So
\begin{equation}
\label{rringgs}
\cO(V)=\cO(Y)^G=(\cO(Z_1) \times...\times \cO(Z_n))^G=(\cO(Z)^{G'}\times...\times\cO(Z)^{G'})^{G/G'}\end{equation}
where $\cO(Z_i)^{G'}\simeq \cO(Z)^{G'}$ via {\it any} $\gamma\in G$ such that $\gamma Z=Z_i$
and $G/G'$ acts on the product via the corresponding permutation representation. Since the last ring in (\ref{rringgs}) contains $\cO(Z)^{G'}$ embedded diagonally it follows that $\cO(Z)^{G'}=\cO(V)$.
\qed

Here is the main result of this section.

\begin{theorem}
\label{main}
Let $X=Spec\ S$ be an ordinary framed curve. Then
the ring $\overline{\tS^{\infty}}$ is a quotient of an ind-\'{e}tale  $\bZ_p^{\times}$-extension of $\overline{S^{\infty}}$.
\end{theorem}

Recall the ring $\stuff$; cf. (\ref{doitild}). By Proposition \ref{escr} we get:

\begin{corollary}
\label{spirt}
\

1) $\overline{\tS^{\infty}}$ is an integral extension of $\overline{S^{\infty}}$.

2) $\stuff$ is an ind-\'{e}tale $\Gamma'$-extension of $\overline{S^{\infty}}$, where $\Gamma'$ is a closed subgroup of $\Gamma:=\bZ_p^{\times}$.
\end{corollary}

For a refinement of this result in the setting of modular curves see Theorem \ref{babanovac}.

For the proof of  Theorem \ref{main}  we need a series of Lemmas.
For the first two Lemmas we let $A$ be a $\d$-ring and we consider the prolongation sequence
$B^r=A[z,z^{-1},z',...,z^{(r)}]\h$. We then denote by $O(r)$ any element of $B^r$.

\begin{lemma}
\label{xlaphi}
Let $\varphi \in A$. Then, for any $n \geq 1$, we have
$$\d^n \left( \frac{z^{\phi}}{z}-\varphi\right)=
z^{-p^n} (z^{(n)})^p-z^{p^{n+1}-2p^n} z^{(n)} +O(n-1)+pO(n+1).$$
\end{lemma}

{\it Proof}. For $\varphi=0$ this is \cite{book}, Lemma 5.19.
Assume now $\varphi$ arbitrary. One checks by induction that
$$\d^n(z-\varphi)=\d^n z+U+pV,$$
where $U=O(n-1)$, $V=O(n)$. Replacing $z$ by $\frac{z^{\phi}}{z}$ we get
$$
\d^n \left( \frac{z^{\phi}}{z}-\varphi\right) = \d^n \left( \frac{z^{\phi}}{z}\right)+
U\left(\frac{z^{\phi}}{z},...,\d^{n-1}\left( \frac{z^{\phi}}{z} \right)\right)
 +pV\left(\frac{z^{\phi}}{z},...,\d^{n}\left( \frac{z^{\phi}}{z} \right)\right),$$
and we conclude by the case $f=0$ of the Lemma.
\qed

\begin{lemma}
\label{logder}
Let $\lambda = 1+p^n a$, $a\in \bZ$. Then
$$\d^n(\lambda z) = z^{(n)}+az^{p^n}+pO(n).$$
\end{lemma}

{\it Proof}. An easy exercise. See also \cite{book}, p. 79.

It is also convenient to formulate the following:

\begin{lemma}
\label{artin}
Let $\Sigmat$ be a ring of characteristic $p$ and consider the $\Sigmat$-algebra
$\Sigmat':=\Sigmat[u]/(u^p-u-G)$ where $G \in \Sigmat$. Consider the action of $\bZ/p\bZ=\{\overline{a}\ ;\ a=0,...,p-1\}$
on $\Sigmat[u]$ defined by $\overline{a} \cdot u=u+\overline{a}$ and consider the induced $\bZ/p\bZ$-action
on $\Sigmat'$. Then any
$\bZ/p\bZ$-invariant element of $\Sigmat'$ is in $\Sigmat$.
\end{lemma}

{\it Proof}.
Let $c \in \Sigmat'$ be the class of $u$. Then $\Sigmat'$ is a free $\Sigmat$-module with basis $1,c,...,c^{p-1}$.
Assume $\sum_{i=0}^{p-1} \lambda_i c^i\in \Sigmat'$ is $\bZ/p\bZ$-invariant, where $\lambda_i \in \Sigmat$.
We want to show that $\lambda_i=0$ for $i\geq 1$. We may assume $\lambda_0=0$. Assume there is a $s \geq 1$
such that $\lambda_s\neq 0$ and let $s$ be maximal with this property. Then
$$\lambda_s(c+1)^s+\lambda_{s-1}(c+1)^{s-1}+...=\lambda_sc^s+\lambda_{s-1}c^{s-1}+...$$
Picking out the coefficient of $c^{s-1}$ we get $s\lambda_s=0$ hence $\lambda_s=0$, a contradiction.
\qed

\bigskip

{\it Proof of Theorem \ref{main}}. For $r \geq 1$ set
$$N^r:=\frac{M^r}{(f-1,\d(f-1),...,\d^{r-1}(f-1))}.$$
Note that
$$\cE^i(\d^{i-1}(f-1))=
\d^{i-1}(\cE^1(f-1))=\d^{i-1}(0)=0.$$
So there are surjective homomorphisms
$N^r \ra \tS^r$, hence surjective homomorphisms $\overline{N^r}\ra \overline{\tS^r}$, hence a surjective homomorphism
\begin{equation}
\label{kki}
\lim_{\ra} \overline{N^r}\ra \lim_{\ra}\overline{\tS^r}=\overline{\tS^{\infty}}.
\end{equation}
Now let $X=\bigcup_{\alpha} X_{\alpha}$, $X_{\alpha}=Spec\ S_{\alpha}$, be an affine open covering
such that $L$ is trivial on each $X_{\alpha}$.
Let $x_{\alpha}$ be a basis of $L$ on $X_{\alpha}$ and
let $z_{\alpha}=x_{\alpha}^{-1}$. Set
$$\begin{array}{lll}
S^r_{\alpha}&:=&S^r_{X_{\alpha}}=(S^r \otimes_S S_{\alpha})\h\\
\  & \  & \  \\
M^r_{\alpha}&:=&M^r_{X_{\alpha}}=(M^r \otimes_S S_{\alpha})\h\end{array}.$$
Then we have an identification
$$M^r_{\alpha}=S^r_{\alpha}[z_{\alpha},z_{\alpha}^{-1},z'_{\alpha},...,z_{\alpha}^{(r)}]\h.$$
Write $f=\varphi_{\alpha} x_{\alpha}^{\phi-1}$, with $\varphi_{\alpha} \in S_{\alpha}^1$.
Since $f$ and $x_{\alpha}$ are invertible in $M_{\alpha}^1$ it follows that $\varphi_{\alpha}$ is invertible in $M_{\alpha}^1$, hence in $S_{\alpha}^1$. Set $N^r_{\alpha}=(N^r \otimes_S S_{\alpha})\h$; hence
$$N^r_{\alpha}:=
\frac{S^r_{\alpha}[z_{\alpha},z_{\alpha}^{-1},z'_{\alpha},...,z_{\alpha}^{(r)}]\h}{\left( \frac{z_{\alpha}^{\phi}}{z_{\alpha}}-\varphi_{\alpha},\d
\left( \frac{z_{\alpha}^{\phi}}{z_{\alpha}}-\varphi_{\alpha}\right),...,\d^{r-1}\left( \frac{z_{\alpha}^{\phi}}{z_{\alpha}}-\varphi_{\alpha}\right)\right)}.$$
For $i \geq 1$ set $u_{i,\alpha}:=\frac{z_{\alpha}^{(i)}}{z_{\alpha}^{p^i}}$.
Also, for $r \geq 1$, set
\begin{equation}
\label{piseout}
\Sigmat^{r,0}_{\alpha}:=\frac{\overline{S^r_{\alpha}}[z_{\alpha},z_{\alpha}^{-1}]}{(z_{\alpha}^{p-1}-
\overline{\varphi_{\alpha}})}=\frac{\overline{S^r_{\alpha}}[z_{\alpha}]}{(z_{\alpha}^{p-1}-\overline{\varphi_{\alpha}})}.\end{equation}
(The latter equality is true because $\varphi_{\alpha}\in (S^1_{\alpha})^{\times}$.)
Then, by Lemma \ref{xlaphi} we have
$\overline{N^1_{\alpha}}=\Sigmat^{1,0}_{\alpha}[u_{1,\alpha}]$
and
$$\overline{N^r_{\alpha}}=\frac{\Sigmat^{r,0}_{\alpha}[u_{1,\alpha},...,u_{r,\alpha}]}{(u_{1,\alpha}^p-u_{1,\alpha}-G_0,...,
u_{r-1,\alpha}^p-u_{r-1,\alpha}-G_{r-2})},\ \ \ r \geq 2,$$
 where
$G_0\in \Sigmat_{\alpha}^{r,0}$,
and
$$G_i \in \Sigmat_{\alpha}^{r,i}:=\frac{\Sigmat^{r,0}_{\alpha}[u_{1,\alpha},...,u_{i,\alpha}]}{(u_{1,\alpha}^p-u_{1,\alpha}-G_0,...,u^p_{i,\alpha}-u_{i,\alpha}-
G_{i-1})},\ \ \ i \geq 1.$$
Clearly the  schemes $Spec\ \Sigmat^{r,i}_{\alpha}$, for various $\alpha$s naturally glue to give a scheme $Spec\ \Sigmat^{r,i}$; so
$\Sigmat^{r,i}\otimes_{\overline{S}}\overline{S_{\alpha}}=\Sigmat^{r,i}_{\alpha}$ for all $\alpha$.
Note that we have
\begin{equation}
\label{axes}
\Sigmat^{r,i}_{\alpha}=\frac{\Sigmat^{r,i-1}_{\alpha}[u_{i,\alpha}]}{(u_{i,\alpha}^p-u_{i,\alpha}-G_{i-1})}\end{equation}
and
natural inclusions
\begin{equation}
\Sigmat^{r,0}_{\alpha}\subset \Sigmat^{r,1}_{\alpha} \subset ... \subset \Sigmat^{r,r-1}_{\alpha} \subset \overline{N^r_{\alpha}}=\Sigmat^{r,r-1}_{\alpha}[u_{r,\alpha}].\end{equation}
So we have natural homomorphisms
$$...\ra \Sigmat^{r,r-1}_{\alpha} \ra \overline{N^r_{\alpha}} \ra \Sigmat^{r+1,r}_{\alpha}\ra \overline{N^{r+1}_{\alpha}}\ra ...$$
which shows that, for each $\alpha$,
$$(\lim_{\stackrel{\ra}{r}} \overline{N^r})\otimes_{\overline{S}}\overline{S_{\alpha}}=\lim_{\stackrel{\ra}{r}} \overline{N^r_{\alpha}}=\lim_{\stackrel{\ra}{r}} \Sigmat^{r,r-1}_{\alpha}
=(\lim_{\stackrel{\ra}{r}} \Sigmat^{r,r-1})\otimes_{\overline{S}}\overline{S_{\alpha}}.$$
These isomorphisms glue together to give an isomorphism
$$\lim_{\ra} \overline{N^r}=\lim_{\ra} \Sigmat^{r,r-1}.$$
We are left to proving that $\lim_{\ra} \Sigmat^{r,r-1}$ is an ind-\'{e}tale  $\bZ_p^{\times}$-extension of $\overline{S^{\infty}}=\lim_{\ra} \overline{S^r}$.

Start by noting  that the maps $\Sigmat^{r,r-1}_{\alpha}\ra \Sigmat^{r+1,r}_{\alpha}$ are injective.
Also $\overline{S^r_{\alpha}}\ra \Sigmat^{r,r-1}_{\alpha}$ are injective and \'{e}tale; cf. (\ref{piseout}) and (\ref{axes}).
Now the group $\Gamma=\bZ_p^{\times}$ acts on $M^r_{\alpha}$ via the rule
$\gamma \cdot z_{\alpha}^{(i)}=\d^i(\gamma z_{\alpha})$ for $\gamma \in \Gamma$. This induces a $\Gamma$-action on $N^r_{\alpha}$
and hence a $\Gamma$-action on $\overline{N^r_{\alpha}}$.  The latter factors through an action of
$\Gamma_r:=(\bZ/p^{r+1}\bZ)^{\times}$. Moreover, for $i \leq r-1$, $\Sigmat^{r,i}_{\alpha}$ is $\Gamma_r$-stable and the $\Gamma_r$-action on $\Sigmat^{r,i}_{\alpha}$
 factors through a $\Gamma_i$-action. For a fixed  $r$  we will prove by induction on $0 \leq i\leq r-1$ that
\begin{equation}
\label{indoi}
(\Sigmat^{r,i}_{\alpha})^{\Gamma_i}=\overline{S^r_{\alpha}}.\end{equation}
This will end the proof of the Theorem; indeed from the above we trivially get that the maps $\Sigmat^{r,r-1}\ra \Sigmat^{r+1,r}$ are injective,
the maps $\overline{S^r}\ra \Sigmat^{r,r-1}$ are injective and \'{e}tale, and, with respect to the induced action,
$$(\Sigmat^{r,r-1})^{\Gamma_{r-1}}=\overline{S^r},$$
showing that $\lim_{\ra} \Sigmat^{r,r-1}$ is an ind-\'{e}tale  $\bZ_p^{\times}$-extension of $\overline{S^{\infty}}=\lim_{\ra} \overline{S^r}$.

Let us check (\ref{indoi}).
For $i=0$  we proceed as follows. Let $b \in \Sigmat^{r,0}_{\alpha}$ be the class of $z_{\alpha}$ and let
$\Gamma_0=\bF_p^{\times}=\langle \zeta \rangle$, $\zeta$ a primitive root.  Then $\Sigmat^{r,0}_{\alpha}$ is a free $\overline{S^r_{\alpha}}$-module with  basis $1,b,b^2,...,b^{p-2}$.
If $\sum_{l=0}^{p-2}\lambda_l b^l$ is $\Gamma_0$-invariant
 (where $\lambda_l \in \overline{S^r_{\alpha}}$)
  then $\sum_{l=0}^{p-2} \lambda_l \zeta^l b^l=\sum_{l=0}^{p-2}\lambda_l b^l$. Since $\zeta$ is primitive we get
 $\lambda_1=...=\lambda_{p-2}=0$,
 and the case $i=0$ is proved.

Now assume $(\Sigmat^{r,i-1}_{\alpha})^{\Gamma_{i-1}}=\overline{S^r_{\alpha}}$ and let us prove (\ref{indoi}).
Recall the equation \ref{axes} and consider the subgroup
$$\Delta_i:=\{\gamma_0,...,\gamma_p\}\subset \Gamma_i,\ \ \ \gamma_a=1+p^ia+p^{i+1}\bZ;$$
so $\Delta_i$ is isomorphic to $\bZ/p\bZ$ via $\gamma_a \mapsto \overline{a}$.
Note that $\Delta_i$ acts trivially on $\Sigmat^{r,i-1}_{\alpha}$.
By Lemma \ref{logder} the $\Delta_i$-action on $\Sigmat^{r,i}_{\alpha}$ corresponds to the $\bZ/p\bZ$-action
induced by $\overline{a}\cdot u_i=u_i+\overline{a}$, so we are in the situation described in Lemma \ref{artin} and we may conclude by that Lemma plus the equality $(\Sigmat^{r,i-1}_{\alpha})^{\Gamma_{i-1}}=\overline{S^r_{\alpha}}$. This ends the proof of (\ref{indoi}) and hence of the Theorem.
\qed

\begin{remark}
\label{simaiintegral}
Assume that the following conditions are satisfied. (This is the case, as we shall see, in the applications to modular curves.)

1) The element $\overline{\varphi} \in \overline{S^1}$ belongs to $\overline{S^0}$.

2) The polynomial $z^{p-1}-\overline{\varphi}$ is irreducible in $Frac(\overline{S^0})[z]$.

\noindent (Here $Frac$ means {\it field of fractions}.) Then the field
\begin{equation}
\label{igusafield}
\frac{Frac(\overline{S^0})[z]}{(z^{p-1}-\overline{\varphi})}\end{equation}
 is embedded into the ring $\overline{\tS^{\infty}}\otimes_{\overline{S^0}}Frac({\overline{S^0}})$.
\end{remark}

\subsubsection{The rings $\tM^r$} For our applications to Igusa modular forms we need one more general construction.
Given a framed curve $X=Spec\ S$ we define
$$\begin{array}{lll}
\tM^r & := & (\tS^r \otimes_{S^r} M^r)\h\\
\  & \  & \  \\
\tM^{\infty} & := & \lim_{\ra} \tM^r.\end{array}$$

The inclusion $\tS^r \subset S^r_{for}$ and the homomorphism $\cE^r:M^r \ra S^r_{for}$ induce a homomorphism (still denoted by $\cE^r$ and still referred to as {\it $\d$-expansion map}),
$$\cE^r:\tM^r \ra S^r_{for}.$$
Assume now $L$ is trivial on $X$; if $x$ is a basis of $L$ then
$$\tM^r=\tS^r[x,x^{-1},x',...,x^{(r)}]\h.$$
Define, in case $L$ is trivial,
$$\tM^r(w):=\tS^r x^w \subset \tM^r;$$
the latter definition is independent of the choice of the basis $x$.
(We will not need and hence we will not define the space $\tM^r(w)$ in case $L$ is not trivial.)
Note also that the element
\begin{equation}
\label{f0}
f^0:=\cE(x^{-1}) \otimes x=\cE(x^{-1})x\in \tM^0(1),\end{equation}
image of the tautological element $x^{-1}\otimes x \in M\otimes_SM$,
does not depend on the choice of the basis $x$ of $L$ and has the property that
$$\cE^0(f^0)=1.$$
We may refer to $f^0$ as the {\it tautological element} of $\tM^0(1)$.

\begin{lemma}
\label{romas}
For any $g \in M^r(w)$ we have $(f^0)^{-w}g =E^r(g)$ in $\tM^r(0)=\tS^r$.
\end{lemma}

{\it Proof}.
If $g=\gamma x^w$, $\gamma \in S^r$, then the image of $g$ in $\tM^r$ equals $E^r(\gamma)x^w$ hence the image of $(f^0)^{-w}g$ in $\tM^r$ equals
$$(E(x^{-1})x)^{-w}E^r(\gamma)x^w=E^r(x^w)x^{-w}E^r(\gamma)x^w=E^r(\gamma x^w)=E^r(g).$$
\qed

The following is trivial to check:

\begin{proposition}
\label{Kracow}
Assume $L$ is trivial on $X$. Then:

1) The homomorphisms $\overline{M^r}\ra \overline{\tM^r}$,  $\overline{M^{\infty}}\ra \overline{\tM^{\infty}}$
are injective. In particular the homomorphisms $M^r \ra \tM^r$, $M^{\infty}\ra \tM^{\infty}$ are injective with torsion free cokernel.

2) The homomorphisms $\tM^r(w)\ra S^r_{for}$ are injective.

3) The homomorphisms $\tM^r \ra \tM^{r+1}$ are injective.
\end{proposition}

\subsection{Application to $\d$-modular forms}
In what follows we shall apply the previous construction to certain framed curves arising from modular (respectively Shimura) curves.

 \subsubsection{Review of modular curves}
 \label{modu}
 Our reference here is \cite[p.~112]{DI}.
The discussion there
involves the modular curve parameterizing elliptic curves with an
embedding of $\mu_N$ rather than $\bZ/N\bZ$ as here.
But, the two modular curves are isomorphic over $\ZN$: see \cite[p.~113]{DI}.

Let $N>3$ be an integer. Consider the $\Z[1/N,\zeta_N]$-scheme $$Y_{\bZ[1/N,\zeta_N]}:=Y_1(N)_{\bZ[1/N,\zeta_N]}$$
 representing the functor taking
a $\Z[1/N,\zeta_N]$-algebra $B$ to the set of isomorphism classes of pairs
$(E,\alpha)$ where $E$ is an elliptic curve over $B$
and $\alpha$ is a level $\Gamma_1(N)$-structure (i.e. $\alpha\colon (\Z/N\Z)_B \subset E[N]$ is an inclusion).
Inclusions as above are the same as inclusions $\mu_{N,B} \subset E[N]$ because we have fixed an $N$th root of unity.
Denote by $L_{\bZ[1/N,\zeta_N]}$ the invertible
sheaf on $Y_{\bZ[1/N,\zeta_N]}$ obtained by taking the direct image of the sheaf
of relative $1$-forms on the universal elliptic curve over $Y_{\bZ[1/N,\zeta_N]}$.
Furthermore consider the
 $\Z[1/N,\zeta_N]$-scheme $X_1(N)_{\bZ[1/N,\zeta_N]}$ defined by taking  the Deligne-Rapoport compactification of $Y_{\bZ[1/N,\zeta_N]}$:
see~\cite[pp.~78--81]{DI} and still denote by $L_{\bZ[1/N,\zeta_N]}$ the natural extension to $X_1(N)$ with the  property
that for  any $\bZ[1/N,\zeta_N]$-algebra $B$ the $B$-module
$M(B,\kappa,N)$
of {\it modular forms over $B$ of level $\Gamma_1(N)$ and weight $\kappa$}
(in the sense of \cite{DI}, p. 111), identifies with the space of sections
$H^0(X_1(N)_B,L_B^{\otimes \kappa})$,
 where $L_B$
is the corresponding sheaf obtained by pull-back.
In particular the normalized Eisentein forms $E_4,E_6.E_{p-1}$ belong to the spaces
$M(\bZ_p,4,N)$, $M(\bZ_p,6,N)$, $M(\bZ_p,p-1,N)$ respectively.
The cusp $\infty$ on $X_1(N)_{\bZ[1/N,\zeta_N]}$ is a $\ZN$-point
$$P_{\bZ[1/N,\zeta_N]}\in X_1(N)_{\bZ[1/N]}(\bZ[1/N,\zeta_N]).$$
 Furthermore, for $p$ not dividing $N$, we choose a homomorphism
$\ZN \ra R=R_p=\hat{\bZ}^{ur}_p$ and we  denote by $Y_R,L_R,P_R$ (or simply by $Y,L,P$) the objects over $R$ obtained from $Y_{\bZ[1/N]},L_{\bZ[1/N]},P_{\ZN}$
by base change.
We denote by
$$Y_{ord}\subset Y=Y_R=Y_1(N)_R$$
the {\it ordinary locus} of $Y$ i.e. the locus in $Y$ where $E_{p-1}$ is invertible.
Finally we let $X=X_R$ be an arbitrary affine open set of $Y=Y_R$. (We will often assume in what follows that $X \subset Y_{ord}$ but at this stage we do not make this assumption.)

This provides us with data (\ref{data1}), (\ref{data5}), where $X^*=X_1(N)_R$. By \cite{Katz} we also have at our disposal an expansion map
$E:M \ra S_{for}$ as in (\ref{data6}), the {\it Fourier expansion map}; recall from \cite{Katz} that this map is obtained by interpreting sections of powers of $L$ as modular forms and evaluating modular forms on the Tate curve
 $\operatorname{Tate}(q)/R((q))$
equipped with its standard $1-$form and its
 standard immersion of $\mu_{N,R} \simeq (\Z/N\Z)_{R}$. A tuple $X,L,P,E$ as above is a Fourier-framed curve
 and $X$ will be referred to as a {\it modular Fourier-framed curve}.

 Alternatively we may keep the data $X,L$ above as our data (\ref{data1}) but instead of $P$ above we may choose
 a point $P \in X(R)$ represented by an elliptic curve $E$ which is a canonical lift of an ordinary
 elliptic curve $\bar{E}$ over $k$. This provides us with data (\ref{data3}). Fix a  $\Z_p$-basis of the Tate module $T_p(\bar{E})$. This basis
 defines an isomorphism between $Spf\ R[[t]]$ and
the completion of $X$ along  the image of $P$.
The Serre-Tate parameter $q$ corresponds to the value
of $1+t$. As in \cite{mori} (cf. also \cite{book}, p.252), we have at our disposal an expansion map $E:M \ra S_{for}$
(\ref{data6}). A tuple $X,L,P,E$ as above is a Serre-Tate-framed curve
 and $X$ will be referred to as a {\it modular Serre-Tate-framed curve}.

In either case (Fourier and Serre-Tate) we have a natural homomorphism
$$\bigoplus_{\kappa\geq 0} M(R,\kappa,N) \ra M=M_X.$$

 \subsubsection{Review of Shimura curves}
 \label{shimu}
 The references here are \cite{shimura,book,buzzard}.

 Let $D$ be a non-split indefinite quaternion algebra over $\bQ$.
Fix a maximal order $\cO_D$ once and for all.
Let $X^D(\cU)_{\bZ[1/m]}$ be the Shimura curve attached to the pair $(D,\cU)$,
where $\cU$ is a sufficiently small
compact subgroup of $(\cO_D \otimes (\varprojlim \bZ/n\bZ))^{\times}$
such that $X^D(\cU)$ is connected and $m \in \bZ_{>0}$ is an appropriate integer; see~\cite{buzzard}.
If
$D$ and $\cU$ satisfy the conditions in~\cite{buzzard}; then for
some $m$ the Shimura curve $X^D(\cU)_{\bZ[1/m]}$ is a
$\Z[1/m]$-scheme with geometrically integral fibers, such that for
any $\Z[1/m]$-algebra $B$ the set $X^D(\cU)_{\bZ[1/m]}(B)$ is in bijection with the
set of isomorphism
  classes of triples $(E,i,\alpha)$ where
$(E,i)$ is a false elliptic curve over $B$ (i.e. $E/B$ is an
abelian scheme of relative dimension $2$ and $i\colon \cO_D \ra
End(E/B)$ is an injective ring homomorphism)
  and $\alpha$ is a level $\cU$ structure.
  Let $Y_{\bZ[1/m]}$ be an affine open subscheme of $X^D(\cU)_{\bZ[1/m]}$.
  For any $p$ not dividing $m$ denote by $Y=Y_R$  the $R$-scheme obtained by base change from $Y_{\bZ[1/m]}$, where $R=R_p=\hat{\bZ}^{ur}_p$.
    We denote by $L=L_R$ the invertible sheaf on $Y_R$ of {\it false $1$-forms}; cf. \cite{book}, p.230.
Finally let $X=X_R \subset Y_R$ be an arbitrary affine open set.
The pair $X,L$ provides us with the data (\ref{data1}).

Next,
by the proof of Lemma~2.6 in~\cite{shimura}, for any $p$ sufficiently big,
there exist infinitely many $k$-points $\bar{P} \in X^D(\cU)(k)$
whose associated triple
$(\bar{E},\bar{i},\bar{\alpha})$ is such that
\begin{enumerate}
\item $\bar{E}$ is ordinary, and
\item if $\bar{\theta}$ is the unique principal polarization
compatible with $\bar{i}$, then $(\bar{E},\bar{\theta})$
is isomorphic to the polarized Jacobian of a genus-$2$ curve.
\end{enumerate}
Choose a $p$ and a point
$\bar{P} \in X(k)$ as above.
Let $E$ be the canonical lift of $\bar{E}$.
Since $End(E) \simeq End(\bar{E})$,
the embedding $\bar{i}\colon \cO_D \to End(\bar{E})$
induces an embedding $i\colon \cO_D \to End(E)$.
Also the level $\cU$ structure
$\bar{\alpha}$ lifts to a level $\cU$ structure on $(E,i)$.
Let $P:=(E,i,\alpha) \in X(R)$; this provides us with the data (\ref{data3}).
Let $\bar{E}^\vee$ be the dual of $\bar{E}$.
By Lemma 2.5 in \cite{shimura}, there exist
$\Z_p$-bases of the Tate modules $T_p(\bar{E})$ and $T_p(\bar{E}^\vee)$,
corresponding to each other under $\bar{\theta}$,
such that any false elliptic curve over $R$ lifting $(\bar{E},\bar{i})$ has a
diagonal Serre-Tate matrix $\operatorname{diag}(q,q^{disc(D)})$
with respect to these bases.
Fix such bases.
They define an isomorphism between $Spf\ R[[t]]$ and
the completion of $X$ along  the image of $P$.
The Serre-Tate parameter $q$ corresponds to the value
of $1+t$. As in \cite{book}, p.252, we have at our disposal an expansion map
as in (\ref{data6})
(referred to in \cite{book} as a {\it Serre-Tate expansion map}).
A tuple $X,L,P,E$ as above is a Serre-Tate-framed curve
 and $X$ will be referred to as a {\it Shimura  Serre-Tate-framed curve} or simply as a
 {\it Shimura framed curve} (for there is no Fourier side of the story in the Shimura curve case).

\subsubsection{$\d$-modular functions and forms}
\label{seinf}
 In what follows we assume  $X$ is either a modular (Fourier or Serre-Tate) framed curve (cf. section \ref{modu})  or a Shimura framed curve
 (cf.  section \ref{shimu}).

Following \cite{book} the elements of $M^r$ will be called {\it $\d-$modular functions}
(holomorphic on $X$) and the elements of $M^r(w)$ will be called {\it $\d$-modular forms of weight $w$}
(holomorphic on $X$). The $\d$-expansion maps are called $\d$-{\it Fourier} respectively $\d$-{\it Serre-Tate} expansion maps.

\begin{remark}
\

1) Assume $X=Y=Y_1(N)_R$ and consider ``test objects" of the form $(E,\alpha,\omega,S^*)$
 where $S^*=(S^n)$ is a prolongation sequence over $R$ of $p$-adically complete rings,  $E/S^0$ is an elliptic curve, $\alpha$ is a $\Gamma_1(N)$-level structure on $E$, and $\omega$ is an invertible $1$-form on $E$. Then giving an element of $M^r_X$ is the same as giving a rule $f$ that attaches to any test object as above   an element
$$f(E,\alpha,\omega,S^*)\in S^r$$ which only depends on the isomorphism class of $(E,\alpha,\omega)$ and is functorial in $S^*$. Moreover the elements of $M^r_X(w)$ correspond to rules $f$ such that
\begin{equation}
\label{bobs}
f(E,\alpha,\lambda \omega,S^*)=\lambda^{-w}f(E,\alpha,\omega,S^*)
\end{equation}
for $\lambda\in (S^0)^{\times}$. So the space $M^r_X(w)$ identifies with the space $M^r(\Gamma_1(N),R,w)$ in \cite{eigen}.

2) Assume $X=Y_{ord}\subset Y$ is the  locus in $Y=Y_1(N)_R$ where $E_{p-1}$ is invertible. Then giving an element in $M^r_X$ is the same as giving a rule $f$ as above but only defined on test objects $(E,\alpha,\omega,S^*)$ such that $E$ has ordinary reduction mod $p$. Then the elements of $M^r_X(w)$ correspond to the rules (defined for elliptic curves with ordinary reduction) satisfying equation (\ref{bobs}). So the space $M^r_X(w)$ identifies, in this case, with the space $M^r_{ord}(\Gamma_1(N),R,w)$ in \cite{eigen}.
\end{remark}

Going back to the case of a general modular or Shimura framed curve we have, under certain additional assumptions,  some remarkable forms
$f^1,f^{\partial}, f^{\sharp}$,
respectively
$f^0, f^{\natural}$,
which we now review (respectively introduce).

\subsubsection{Review of the forms $f^1,f^{\partial}$}
\label{xxx}
The references here are \cite{difmod,Barcau,shimura,book}.

In the Fourier expansion case we set
\begin{equation}
\label{Pssi}
\Psi:=\frac{1}{p} \log \frac{q^{\phi}}{q^p}
    := \sum_{n \geq 1}
        (-1)^{n-1}n^{-1} p^{n-1} \left( \frac{q'}{q^p} \right)^n
    \in S^1_{for}=R((q))\h[q']\h.\end{equation}
     In the Serre-Tate expansion case we still denote by
    \begin{equation}\label{Pssii}
    \Psi\in S^1_{for}=R[[t]][t']\h\end{equation}
the image of $\Psi$ in (\ref{Pssi}) via the homomorphism
\begin{equation}
\label{micar}
R((q))\h[q']\h \ra R[[t]][t']\h\end{equation}
given by $q\mapsto t+1$, $q'\mapsto \d(t+1)$. So $\Psi$ in (\ref{Pssii}) is given by
\begin{equation}
\label{Pssiii}
\Psi=\frac{1}{p}(\phi-p)\log(1+t).\end{equation}

In the next two Propositions $X=Spec\ S$ is either a modular framed curve or a Shimura framed curve; recall the $\d$-expansion maps $\cE^r:M^r \ra S^r_{for}$, cf. (\ref{adelita}).

\begin{proposition}
\label{funu}
\cite{difmod, shimura} There exists a unique form $f^1 \in M^1(-1-\phi)$ such that
$$\cE^1(f^1)=\Psi.$$
\end{proposition}

\begin{proposition}
\label{fpartz}
\cite{Barcau,shimura,book}
Assume the reduction mod $p$ of $X$, $\overline{X}$, is contained in the ordinary locus of the modular (respectively Shimura) curve.
Then there exists a unique form $f^{\partial} \in M^1(\phi-1)$ which is invertible in the ring $M^1$
such that
$$\cE^1(f^{\partial})=1.$$
Furthermore its reduction mod $p$, $\overline{f^{\partial}}\in \overline{M^1(\phi-1)}$, coincides with the image of the Hasse invariant  $\overline{H} \in \overline{M^0(p-1)}$. In particular $X$ is an ordinary framed curve in the sense of Definition \ref{ass2}.
\end{proposition}

So Theorem \ref{main} and Corollary \ref{spirt} hold for any $X$ as in Proposition \ref{fpartz}.
For the definition of the Hasse invariant in the modular curve case see \cite{Katz}.
For the definition in the Shimura curve case see \cite{book}, p. 236.
Proposition \ref{fpartz} justifies the use of the  terminology ``ordinary" introduced in Definition \ref{ass2}.
The condition that
 $\overline{X}$ be contained in the ordinary locus of the modular (respectively Shimura) curve is satisfied if we assume that  a lift $H \in M^0(p-1)$ of the Hasse invariant
is invertible on $X$; in the case of modular curves (cf. section \ref{modu}) this is so if we assume that the Eisenstein
form $E_{p-1} \in M(R,p-1,N)$ is invertible on $X$, i.e. if $X \subset Y_{ord}$.

\subsubsection{Igusa $\d$-modular functions and forms}
\label{pa}
  Let $X$ be  either  modular framed curve (cf. section \ref{modu})  or a Shimura framed curve (cf.  section \ref{shimu}).
Assume
 $\overline{X}$ is contained in the ordinary locus of the modular (respectively Shimura) curve and assume $L$ is trivial on $X$.
The elements
of $\tM^r$ (respectively $\tM^r(w)$) will be called {\it Igusa $\d$-modular functions} (respectively
{\it Igusa $\d$-modular forms of weight $w$}). The ``Igusa'' terminology is justified by the fact that,
due to Proposition \ref{fpartz}, the assumptions in Remark \ref{simaiintegral} are satisfied in the modular curve case (of section \ref{modu}) and, in this case, the field  (\ref{igusafield}) identifies with the field
of rational functions on the Igusa curve \cite{gross}, pp. 460-461.

\subsubsection{The forms $f^0$ and $f^{\natural}$} Assume $L$ is trivial on $X$.
We may consider  the Igusa $\d$-modular form
$f^0 \in \tM^0(1)$, cf. (\ref{f0}). So, by Proposition \ref{Kracow},  $f^{\partial}=(f^0)^{\phi-1}$ in $\tM^1$.

\begin{corollary}
\label{partdeux}
$(\overline{f^0})^{p-1}=\overline{H}$ in $\overline{\tM^1}$.
\end{corollary}

{\it Proof}. We have
$$(\overline{f^0})^{p-1}=\overline{(f^0)^{\phi-1}}=\overline{f^{\partial}}=\overline{H}.$$
\qed

Also let us note that we have, at our disposal, a remarkable Igusa $\d$-modular form of weight $0$:
\begin{equation}
\label{natural}
f^{\natural}:=(f^0)^{\phi+1}f^1 \in \tM^1(0)=\tS^1.\end{equation}
This form will play a key role in the several primes theory; cf. (\ref{fimportant}) below.

\subsubsection{Review of the forms $f^{\sharp}$}

The references here are \cite{eigen, local}.

First recall some standard definitions in the classical theory of complex modular forms \cite{DI}.
Consider the space
$$S_m(\Gamma_1(N),\bC)\subset M(\bC,m,N)$$
 of (classical)
  cusp forms $\sum a_n q^n$ of weight $m$ on $\Gamma_1(N)$ over the complex field $\bC$.
   On this space one has
Hecke operators $T_m(n)$ acting, $n \geq 1$. An {\it eigenform}  is a nonzero
simultaneous eigenvector for all $T_m(n)$, $n \geq 1$. An
eigenform $f=\sum_{n \geq 1} a_n q^n$ is {\it
normalized} if $a_1=1$; in this case $T_m(n)f=a_n  f$ for all
$n \geq 1$. One associates to any  eigenform $f$ its {\it system of eigenvalues} $l \mapsto
a_l$, $(l,N)=1$. A {\it newform} is a normalized eigenform whose
 system of eigenvalues does not come from
 a system of eigenvalues associated
to an eigenform of level $M$ with $M\ |\ N$, $M \neq
N$.
We denote by  $S_m(\Gamma_0(N),\bC)$ the subspace
of $S_m(\Gamma_1(N),\bC)$ of all cusp forms of weight $m$ on
$\Gamma_0(N)$. By Eichler-Shimura theory to each newform
$f =\sum a_n q^n \in S_2(\Gamma_0(N),\bC)$
with $\bZ$-Fourier coefficients there corresponds an elliptic curve $E_f$ over $\bQ$.
We say that $f$ is of CM type if $E_f$ has CM.

Now we place ourselves in the context of a modular Fourier-framed curve.

\begin{theorem}
\label{fur}
\cite{eigen}
 Let $f =\sum a_n q^n \in S_2(\Gamma_0(N),\bC)$ be a newform
with $\bZ$-Fourier coefficients. Assume that
 $f$ is not of CM type and $p$ is sufficiently big.
 Then
there is a unique form $f^{\sharp} \in
M^2(0)$  with $\d-$Fourier expansion:
\begin{equation} \label{h1}
\cE^2(f^{\sharp}) = \frac{1}{p} \sum_{n\geq 1}
\frac{a_n}{n} (q^{n\phi^2} -a_p q^{n\phi} +p q^n) \in
 R((q))\hat{\ }[q',q'']\hat{\ }.
 \end{equation}
\end{theorem}

There is a corresponding result for the case
when $f$ is of CM type; cf. \cite{eigen}.

\subsubsection{Canonical derivations on $\tS^{\infty}$}
In  this section we let $X=Spec\ S$ be either a modular framed curve or a Shimura framed curve
and we assume
 $\overline{X}$ is contained in the ordinary locus; we show that $\tS^{\infty}$ (and other rings related to it) carry a complete system of conjugates of a classical derivation operator on modular forms.
Recall the {\it Serre derivation operator}
$\partial: M\ra M$
introduced by Serre and Katz (cf. \cite{Katz} for the modular curve case and also
\cite{book}, p. 254, for the modular/Shimura case). Recall that $M=\bigoplus_{n \in \bZ}L^{\otimes n}$ and $\partial(L^{\otimes n}) \subset L^{\otimes (n+2)}$.  Similarly one can consider
the {\it Euler derivation operator}
$\cD: M \ra M$
defined by requiring that its restriction to $L^{\otimes n}$ is multiplication by $n$.  Finally recall the Ramanujan form $P \in M^0(2)$; cf. \cite{Katz} in the modular case and \cite{book}, Definition 8.33, in the modular/Shimura curve case. (We stress the fact that
$P$ is well defined because
 we are  assuming that  $\overline{X}$ is contained in the ordinary locus of the modular/Shimura curve.) Consider the derivation
$$\partial^*:=\partial+P\cD:M^0 \ra M^0,$$
  and its complete sequence of conjugates $\partial^*_j:M^{\infty}\ra M^{\infty}$, $j \geq 0$.
\begin{lemma}
\label{trav}
The kernel of the $\d$-expansion map $E^{\infty}:M^{\infty}\ra S^{\infty}_{for}$ is sent into itself by  each of the derivations  $\partial^*_j$.
\end{lemma}

{\it Proof}.
Assume first we are in the Serre-Tate case.
It follows  from \cite{book}, Proposition 8.42, that there exist derivations $D_j$
on $S^{\infty}_{for}$  such that the difference
$E^{\infty}\circ \partial_j^*-  D_j \circ E^{\infty}$
vanishes on $\bigoplus_{w \in W(r)}M^r(w)$ for all $r$.
It follows easily that the difference above vanishes on the whole of $M^{\infty}$. But this implies the assertion of the Lemma.
A similar argument holds in the  Fourier case; one uses \cite{Barcau}, Proposition 4.2, instead.
\qed

Taking the derivations (still denoted by $\partial^*_j$) on $\tS^{\infty}$ induced by $\partial^*_j:M^{\infty}\ra M^{\infty}$ (cf. Lemma
\ref{trav}) we get:

\begin{corollary}
\label{evitt}
 There is a (necessarily unique) complete sequence of conjugates
$$\partial_j^*:\tS^{\infty}\ra \tS^{\infty}$$ of $\partial^*:M^0\ra M^0$.
Moreover $\partial^*_j \tS^s \subset \tS^s$ for all $j,s\geq 0$ and $\partial^*_j \tS^s=0$ for $s<j$
\end{corollary}

Uniqueness in the above statement follows from the fact that $\tS^{\infty}$ is (clearly) $p$-adically separated
and topologically $\d$-generated by $M^0$.

One can ask what is the effect of the derivations $\partial^*_j$ when applied to the form
$f^{\natural}$ in (\ref{natural}). Recall that this form was defined only in case $L$ is trivial on $X$.
Clearly $\partial_j^* f^{\natural}=0$ for $j \geq 2$. For $\partial^*_0,\partial^*_1$ we have:

\begin{proposition}
\label{zece}
$$\begin{array}{rcl}
\partial^*_0 f^{\natural} & = & -1\\
\  & \  & \  \\
\partial_1^* f^{\natural} & = &1.\end{array}$$
\end{proposition}

{\it Proof}.
Recall from \cite{book}, pp.271-272, that
$$\begin{array}{rcl}
\partial^*_0 f^1 & = & -\frac{1}{f^{\partial}}\\
\  & \  & \  \\
\partial^*_1 f^1 & = & f^{\partial}\end{array}$$
in $M^{\infty}$. By Lemma \ref{romas} it follows that
$$\partial^*_0 f^{\natural}=\partial^*_0(E^1(f^1))=E^1(\partial^*_0 f^1)=
E^1(-\frac{1}{f^{\partial}})=-1.$$
The computation for $\partial^*_1f^{\natural}$ is similar.
\qed

Now the derivations $\partial_j^*$ in Corollary \ref{evitt} induce $k$-derivations $\overline{\partial_j^*}:\overline{\tS^{\infty}}\ra \overline{\tS^{\infty}}$.

\begin{lemma}
The kernel of the natural map $\overline{E^{\infty}}:\overline{\tS^{\infty}}\ra \overline{S^{\infty}_{for}}$ is sent into itself by each of the derivations
$\overline{\partial_j^*}:\overline{\tS^{\infty}}\ra \overline{\tS^{\infty}}$.
\end{lemma}

{\it Proof}.
Let $D_j$ be as in the proof of Lemma \ref{trav} and let $\overline{D_j}:\overline{S^{\infty}_{for}}\ra \overline{S^{\infty}_{for}}$ be the induced derivations.
So $\overline{E^{\infty}}\circ \overline{\partial_j^*}=\overline{D_j}\circ \overline{E^{\infty}}$ on $\overline{M^{\infty}}$ hence on $\overline{\tS^{\infty}}$.
This implies the statement of the Lemma.
\qed

Recall the ring $\stuff$; cf. (\ref{doitild}).Then we have:

\begin{corollary}
\label{riri}
The derivations $\overline{\partial_j^*}:\overline{\tS^{\infty}}\ra \overline{\tS^{\infty}}$ induce derivations
$\widetilde{\partial_j^*}:\stuff \ra \stuff$.
\end{corollary}

We want to further analyze the derivations in Corollary \ref{riri}. It is convenient to give the following

\begin{definition}
Let $F$ be an algebra over a ring $\Lambda$. (In the applications $\Lambda$ will be either $R$ or $k=R/pR$.) We say that a sequence $\theta_0,\theta_1,\theta_2,...\in Der_{\Lambda}(F,F)$ is
{\it locally finite} if for any $x \in F$ there exists an integer $j_0$ such that for any $j \geq j_0$ we have $\theta_j x=0$.
Given such a locally finite sequence one has an $F$-linear map
\begin{equation}
\label{mapu}
F \times F \times F \times ... \ra Der_{\Lambda}(F,F)\end{equation}
sending any vector $(c_0,c_1,c_2,...)$ into the well defined derivation $\sum_{j \geq 0} c_j \theta_j$.
A  sequence of derivations is called
a {\it pro-basis} of $Der_{\Lambda}(F,F)$ if it is locally finite and the map (\ref{mapu}) is a bijection.
\end{definition}

A basic example is the polynomial ring $F=\Lambda[T,T',T'',...]$; in this case the sequence of derivations
\begin{equation}
\label{lod}
\frac{\partial}{\partial T}, \frac{\partial}{\partial T'}, \frac{\partial}{\partial T''},... \in Der_{\Lambda}(F,F)\end{equation}
is a pro-basis of $Der_{\Lambda}(F,F)$.

Then we have the following:

\begin{proposition}
\label{noua}
Assume $X$ is sufficiently small. (If $X$ is a modular Fourier-framed curve it is enough to assume that
the Eisenstein forms $E_4, E_6, E_{p-1}$ are invertible in $M$.)
 Then the induced derivations
$$\widetilde{\partial_0^*}, \widetilde{\partial_1^*}, \widetilde{\partial_2^*},... \in Der_k(\stuff,\stuff)$$
form a pro-basis of $Der_k(\stuff,\stuff)$.
\end{proposition}

{\it Proof}.
Shrinking $X$ we may assume there is an
 \'{e}tale coordinate $T$ on  $X$. (In the modular curve case, if $E_4,E_6$ are invertible in $M$ then one can take $T=j$, the $j$-invariant.)
  By Corollary
\ref{spirt} and  \cite{book}, Proposition 3.13, the derivations
 (\ref{lod}), for $\Lambda=k$, uniquely  extend to derivations (still denoted by)
$$
\frac{\partial}{\partial T}, \frac{\partial}{\partial T'}, \frac{\partial}{\partial T''},... \in Der_k(\stuff,\stuff).$$
The latter sequence is trivially checked to be a pro-basis of $Der_k(\stuff,\stuff)$. Now express
$$\widetilde{\partial^*_j}=\sum_{s \geq 0} c_{js}\frac{\partial}{\partial T^{(s)}} \in Der_k(\stuff,\stuff),$$
with $c_{js}\in \stuff$. Recall that $\partial^*_j T^{(s)}=0$ for $s<j$. Hence $c_{js}=0$ for $s<j$. Also
$$\widetilde{\partial^*_j}T^{(j)}=\overline{\phi^j(\partial^* T)}=(\overline{\partial^* T})^{p^j} \in \stuff.$$
Note that the image of $\overline{\partial^* T}$ in $\overline{S^{\infty}_{for}}$ is invertible; in the Serre-Tate expansion case this follows from
\cite{book}, Propositions 8.34 and 8.42, while in the case of Fourier expansions this follows from the fact that
$$\partial j\in \bZ_p^{\times} \cdot \frac{E_4^2E_6}{\Delta}.$$
By shrinking $X$ we may assume
  $\overline{\partial^* T}$ is invertible in $\stuff$.
(Again this is automatic in the modular case if $E_4,E_6$ are invertible in $M$.)
In particular $c_{jj} \in (\stuff)^{\times}$ for all $j\geq 0$.
Now for each $j \geq 0$ one can find inductively elements $\gamma_{js}\in \stuff$ with $s \geq j$ such that
$$\frac{\partial}{\partial T^{(j)}}=\sum_{s \geq j} \gamma_{js} \widetilde{\partial_s^*} \in Der_k(\stuff,\stuff).$$
So any derivation
$$\theta=\sum_{j \geq 0} c_j \frac{\partial}{\partial T^{(j)}}\in Der_k(\stuff,\stuff)$$
with $c_j \in \stuff$ can be written as
$$\theta=\sum_{s \geq 0} \left( \sum_{j=0}^s c_j \gamma_{js}\right) \widetilde{\partial^*_s}.$$
This representation is easily seen to be unique  by evaluating successively at $T$, $T'$, $T'$,.... This ends the proof of the Proposition.
\qed

Our next purpose is to explain the relation between our Igusa $\d$-modular functions and Katz's generalized $p$-adic modular functions
\cite{Katzgen, gouvea}.
This relation will play a role later, in the proof of Theorem \ref{babanovac} below.

We first need to review some of Katz's concepts.

\subsubsection{Review of Katz generalized $p$-adic modular functions} The references here are \cite{Katzgen, gouvea}.

Let $B$ be a $p$-adically complete ring, $p\geq 5$, and let $N$ be an integer coprime to $p$.
Consider the functor
\begin{equation}
\label{hite}
\{\text{$p$-adically complete $B$-algebras}\} \ra \{\text{sets}\}
\end{equation}
that attaches to any $A$ the set of isomorphism classes of triples $(E/A,\varphi, \iota)$, where $E$ is an elliptic curve over $A$,
$\varphi$ is a trivialization, and $\iota$ is an arithmetic level $N$ structure. Recall that a {\it trivialization} is an isomorphism
between the formal group of $E$ and the formal group of the multiplicative group; an {\it arithmetic level $N$ structure}
is defined as an inclusion of flat group schemes over $B$ of $\mu_{N}$ into $E[N]$. The functor (\ref{hite})
is representable by a $p$-adically complete ring $\bW(B,N)$.
The elements of this ring are called by Katz \cite{Katzgen} {\it generalized $p$-adic modular forms}. Note that
 $\bW(B,N)=\bW(\bZ_p,N) \widehat{\otimes} B$.
Moreover there is a $\bZ_p^{\times}$-action on $\bW(B,N)$ coming from the action of $\bZ_p^{\times}$ on the formal group of the multiplicative group.
There is a natural {\it Fourier expansion map}
$E:\bW(B,N) \ra B((q))\hat{\ }$ which is injective and has a flat cokernel over $B$. Also $\bW(\bZ_p,N)$ possesses a natural ring endomorphism
$Frob$ which reduces modulo $p$ to the $p$-power Frobenius endomorphism of $\bW(\bZ_p,N) \otimes \bZ/p\bZ$. So if $R=\hat{\bZ}_p^{ur}$, as usual,
and if $\phi$ is the automorphism of $R$ lifting Frobenius then
$Frob \widehat{\otimes} \phi$ is a lift of Frobenius on
$$\bW:=\bW(R,N)=\bW(\bZ_p,N) \widehat{\otimes} R$$ which we denote  by $\phi_0$.
 Moreover the homomorphism $\bW(R,N) \ra R((q))\hat{\ }$ commutes with the action of $\phi_0$ where $\phi_0$  on $R((q))\hat{\ }$
is defined by $\phi_0(\sum a_nq^n):= \sum \phi(a_n)q^{np}$.

For any $\bZ[1/N,\zeta_N]$-algebra $B$  the space $M(B,\kappa,N)$ of modular forms over $B$ of weight $\kappa$ and level $\Gamma_1(N)$ has an embedding
$$M(B,\kappa,N) \subset \bW(B,N).$$
The space $M(B,\kappa,N)$ is stable under the $\bZ_p^{\times}$-action on $\bW(B,N)$ and $\lambda \in \bZ_p^{\times}$ acts on $M(B,\kappa,N)$ via multiplication by $\lambda^k$.
Recall that we denoted by $Y_{ord}\subset Y=Y_R=Y_1(N)_R$
  the locus in $Y$ where the Eisenstein form $E_{p-1}\in M(\bZ_p,p-1,N)$ is invertible. Then, since $E_{p-1}$ is invertible in $\bW$ we get a homomorphism
$$M_{Y_{ord}}=\bigoplus_{k \in \bZ} L_{Y_{ord}}^{\otimes k} \ra \bW.$$
More generally, if $X$ is any affine open subset of $Y_{ord}$ then one can find $g \in M_{Y_{ord}}$ of weight $0$, $\overline{g}\neq 0$,
and a homomorphism
\begin{equation}
\label{mafura}
M:=M_X:=\bigoplus_{\kappa \in \bZ} L_X^{\otimes \kappa} \ra \bW_g=\bW[1/g].
\end{equation}
(So if $X=Y_{ord}$ we may take $g=1$.)
Since $g$ has weight $0$, $\widehat{\bW_g}$ has an induced $\bZ_p^{\times}$-action and
the homomorphism (\ref{mafura}) is $\bZ_p^{\times}$-equivariant if $\lambda\in \bZ_p^{\times}$ acts on each $L^{\otimes k}$ via multiplication by $\lambda^k$.

Finally recall Katz's {\it ring of  divided congruences} \cite{gouvea},
$$\bD:=\bD(R,N):=\{f \in \bigoplus_{\kappa\geq 0} M(R,\kappa,N)\otimes_R K;\ E(f)\in R[[q]]\},$$
where $K:=R[1/p]$.
This ring naturally embeds into Katz's  {\it ring of  holomorphic generalized $p$-adic modular forms},
$$\bV:=\bV(R,N)=\{f \in \bW(R,N);\ E(f)\in R[[q]]\},$$
and the image of $\bD$ in $\bV$ is $p$-adically dense. For simplicity we sometimes identify $\bD,\bV,\bW$ with subrings of $R((q))\h$; i.e.
we view $$\bD \subset \bV \subset \bW \subset R((q))\h.$$

We will need the following:

\begin{lemma} \label{denn}
$\bD+R[\Delta^{-1}]$ is $p$-adically dense in $\bW$.
\end{lemma}

{\it Proof}.
It is enough to check that $\bV+R[\Delta^{-1}]$ is $p$-adically dense in $\bW$.

 We first claim that for any $f \in \bW$ there exists
a sequence of polynomials $F_n\in R[t]$ in a variable $t$ such that $F_{n+1}-F_n \in p^n R[t]$ for $n \geq 0$ and such that
$$f-F_n(\Delta^{-1})\in p^n R((q))\h+R[[q]].$$
To check the claim we construct $F_n$ by induction. We may take $F_0=0$. Now, assuming $F_n$ was constructed, write
$$f-F_n(\Delta^{-1})=p^n G+p^{n+1}H+S,\ \ G \in \sum_{i=1}^N Rq^{-i},\ \ H \in R((q))\h,\ \ S \in R[[q]].$$
Since $\Delta^{-1}-q^{-1}\in R[[q]]$ we can find a polynomial $\Gamma \in R[t]$ of degree $\leq N$ such that
$G-\Gamma(\Delta^{-1})\in R[[q]]$. Then set $F_{n+1}:=F_n+p^n \Gamma$ which ends the inductive step of our construction.

Now let $F_n$ be as in our claim above and set $F:=\lim F_n \in R[t]\h$. Then clearly $f -F(\Delta^{-1}) \in R[[q]] \cap \bW=\bV$.
This implies that $\bV+R[\Delta^{-1}]$ is $p$-adically dense in $\bW$ and we are done.
\qed

\subsubsection{Link between Igusa $\d$-modular functions and Katz generalized $p$-adic modular functions}
Since $\bW:=\bW(R,N)$ is flat over $R$ and comes with a lift of Frobenius,
$\bW$ has a natural structure of $\d$-ring. The same is true about $R((q))\hat{\ }$ and,
moreover the Fourier expansion map $\bW \ra R((q))\hat{\ }$ is a $\d$-ring homomorphism.
Note that $\d$ acts on $R((q))\hat{\ }$ via the operator
\begin{equation}
\label{seface}
\d_0(\sum a_n q^n):=\frac{(\sum \phi(a_n)q^{np})-(\sum a_n q^n)^p}{p};\end{equation}
this is {\it not the same} as the restriction of $\d:S^{\infty}_{for}\ra S^{\infty}_{for}$ to $R((q))\hat{\ }$.
More generally, $\widehat{\bW_g}$  has a structure of $\d$-ring and the induced homomorphism
\begin{equation}
\label{albb}
\widehat{\bW_g} \ra R((q))\hat{\ }\end{equation}
 is a $\d$-ring homomorphism.
(The latter is well defined because the image of $g$ in $k((q))$ is non-zero, hence $g$ is invertible in $R((q))\hat{\ }$.)
Now the reduction mod $p$ of (\ref{albb}) is injective (because the cokernel of $\bW \ra R((q))\hat{\ }$ is flat over $R$) so (\ref{albb}) itself is injective.
We have an induced commutative diagram of $\d$-rings
$$\begin{array}{ccc}
M^{\infty} & \ra & S^{\infty}_{for}\\
\downarrow & \  & \downarrow\\
\widehat{\bW_g} & \ra & R((q))\hat{\ }\end{array}$$
where the right vertical arrow is induced by $q'\mapsto 0$, $q''\mapsto 0$, etc.
Note that $\bZ_p^{\times}$ naturally acts on $M^{\infty}$ compatibly with $\d$ (with $\lambda\in \bZ_p^{\times}$ acting on $L^{\otimes n}$ via multiplication
by $\lambda^n$) and the map $M^{\infty}\ra \widehat{\bW_g}$ is $\bZ_p^{\times}$-equivariant.
Now the kernel of the upper horizontal map
is mapped  to $0$ in $R((q))\hat{\ }$, so is mapped into $0$ in $\widehat{\bW_g}$ (because the lower horizontal arrow is injective).
Hence we have an induced homomorphism
\begin{equation}
\label{kkkk}
\tS^{\infty} \ra \widehat{\bW_g}.
\end{equation}
By Proposition \ref{funu} the kernel of (\ref{kkkk}) contains the images of the elements
\begin{equation}
\label{immages}
f^1,\ \d f^1,\ \d^2 f^1,...\end{equation}
If $L$ is trivial with basis $x$ we set
$$f^{\flat}:=f^1 x^{\phi+1} \in S^1=M^1(0)$$
and then the kernel of (\ref{kkkk}) contains the images of the elements
\begin{equation}
\label{culllpr}
f^{\flat},\d f^{\flat}, \d^2 f^{\flat},...
\end{equation}

In a similar vein
note that
 the kernel of the map
 \begin{equation}
 \label{oldd}
 M^{\infty}\ra \tS^{\infty} \subset S^{\infty}_{for}\end{equation}
  contains
  \begin{equation}
  \label{duss}
f^{\partial}-1,\ \d(f^{\partial}-1),\ \d^2(f^{\partial}-1),....\end{equation}
In the next section we will show, under an appropriate hypothesis, that the kernels of
(\ref{kkkk}) and (\ref{oldd}) are topologically generated by the elements (\ref{culllpr}) and (\ref{duss}) respectively.

\subsubsection{Refinement of results in the modular case}
Using the link between Igusa $\d$-modular functions and Katz generalized $p$-adic modular forms we can refine Theorem \ref{main} in the modular case as follows:

\begin{theorem}
\label{babanovac}
Assume $X=Spec\ S$ is a modular Fourier-framed curve with $E_{p-1}$ invertible on $X$.  The following hold:

1) The map $\overline{\tS^{\infty}}\ra \overline{S^{\infty}_{for}}$ is injective; in particular $\overline{\tS^{\infty}}$ is an integral domain,
and
the map $\overline{\tS^{\infty}}\ra \stuff$ is an isomorphism. Moreover
the ring $\overline{\tS^{\infty}}$ is an ind-\'{e}tale $\bZ_p^{\times}$-extension of $\overline{S^{\infty}}$.

2) The kernel of $\overline{M^{\infty}} \ra \overline{\tS^{\infty}}$ is generated by
$$\overline{f^{\partial}-1},\ \overline{\d(f^{\partial}-1)},\ \overline{\d^2(f^{\partial}-1)},...$$

3) The kernel of $\overline{\tS^{\infty}}\ra \overline{\bW}_g$ is generated by the images of
$$\overline{f^1},\ \overline{\d f^1}, \ \overline{\d^2 f^1},...$$

4) The kernel of $\overline{M^{\infty}} \ra \overline{\bW}_g$ is generated by the elements
$$\overline{f^{\partial}-1},\ \overline{f^1},\ \overline{\d(f^{\partial}-1)}, \ \overline{\d f^1},\ \overline{\d^2(f^{\partial}-1)}, \ \overline{\d^2 f^1},...$$
\end{theorem}

\begin{corollary}
\label{noname}
Assume $X=Spec\ S$ is a modular Fourier-framed curve  with $E_{p-1}$ invertible on $X$. The following hold:

1) The inclusion $\tS^{\infty}\subset S^{\infty}_{for}$ has torsion free cokernel.

2) The kernel of $M^{\infty}\ra S^{\infty}_{for}$ is the $p$-adic closure of the ideal generated by the elements
$$f^{\partial}-1,\ \d(f^{\partial}-1),\ \d^2(f^{\partial}-1),...$$

3) The kernel of $\tS^{\infty}\ra R((q))\hat{\ }$ is the $p$-adic closure of the ideal generated by the images of the elements
$$f^1,\  \d f^1, \ \d^2 f^1,...$$

4) The kernel of $M^{\infty}\ra R((q))\hat{\ }$ is the $p$-adic closure of the ideal generated by the elements
$$f^{\partial}-1,\ f^1,\ \d(f^{\partial}-1),\ \d f^1, \ \d^2(f^{\partial}-1),\ \d^2 f^1,...$$
\end{corollary}

\begin{remark}
Conclusion  1 in Corollary \ref{noname}  should be viewed as a {\it  $\d$-expansion principle}. Conclusions 2 and 4
 should be viewed as  $\d$-analogues of the theorem of Swinnerton-Dyer and Serre according to which
the kernel of the Fourier expansion map
$$\bigoplus_{\kappa \geq 0} M(\bF_p,\kappa,N)\ra \bF_p[[q]]$$
 is generated by $E_{p-1}-1$; cf. \cite{gross}, p. 459.
\end{remark}

To prove Theorem \ref{babanovac} we need the following:

\begin{lemma}
\label{becar}
Any  ordinary $k$-point of $Y=Y_1(N)_R$ has an affine open neighborhood $X\subset Y$ such that $E_{p-1}$ is invertible on $X$, $L$ is trivial on $X$ and
 the natural homomorphism
$$\overline{S} \ra \frac{\overline{S^r}}{(\overline{f^{\flat}},\overline{\d f^{\flat}},...,\overline{\d^{r-1}f^{\flat}})}$$
is an isomorphism.
\end{lemma}

{\it Proof}.
Same argument as in \cite{local}, Lemma 4.66 (where the special case $r=2$ was considered).
\qed

\medskip

{\it Proof of Theorem \ref{babanovac}.}
We are going to use the notation in the proof of Theorem \ref{main}.
In particular recall the rings $\Sigmat^{r,r-1}$ which are finite \'{e}tale extensions of $\overline{S^r}$, with
\begin{equation}
\label{zty}
(\Sigmat^{r,r-1})^{\Gamma_{r-1}}=\overline{S^r}.
\end{equation}

Note that assertion 4  follows from assertions 2 and 3.

We claim that in order to prove assertions 1 and 2 it is enough to show that all the rings $\Sigmat^{r,r-1}$ are integral domains.
 Indeed if this is so then
$$\Sigmat^{\infty}:=\lim_{\ra} \Sigmat^{r,r-1}$$
is an integral domain. We have surjections
\begin{equation}
\label{surjj}
\Sigmat^{\infty} \ra \overline{\tS^{\infty}}\ra \stuff,\end{equation}
where the last ring is an integral domain. Let $I$ be the kernel of the composition (\ref{surjj}). Since the composition $\overline{S^{\infty}}\ra \Sigmat^{\infty} \ra
\stuff$ is injective (cf. Proposition \ref{expansionprinciple}), upon viewing $\overline{S^{\infty}}$ as a subring of $\Sigmat^{\infty}$,
 it follows that $I \cap \overline{S^{\infty}}=0$.
Since $\Sigmat^{\infty}$ is an integral domain and an integral extension of $\overline{S^{\infty}}$ it follows that $I=0$. This forces the surjections in (\ref{surjj})
to be isomorphisms, and so assertions 1 and 2 of the Theorem  follow.

Next note that since $Spec\ \Sigmat^{r,r-1}$ is \'{e}tale and finite over $Spec\
\overline{S^r}$
and since the latter is smooth over $k$, it follows that $Spec\ \Sigmat^{r,r-1}$ is smooth over $k$ so, in particular its connected components are irreducible
and they are finite and \'{e}tale over $Spec\ \overline{S^r}$. So in order to prove that $\Sigmat^{r,r-1}$
is an integral domain it is enough to prove that $Spec\ \Sigmat^{r,r-1}$ is connected.

Consequently in order to prove the Theorem we need to prove connectivity of $Spec\ \Sigmat^{r,r-1}$ and assertion 3.
We will prove these two facts simultaneously. To prove either of these facts
it is enough to prove that these facts hold for each
 of the open sets
of a given open cover of $X$. So we may assume, after shrinking $X$, that  the conclusion of Lemma \ref{becar} holds for $X$, in particular
$L$ is trivial on the whole of $X$
so $f^{\flat}$ is defined and $f^1$ and $f^{\flat}$ differ by a unit.
Consider the scheme $Spec\ T^r$
defined by the cartesian diagram
$$\begin{array}{rcl}
Spec\ T^r & \ra & Spec\ \Sigmat^{r,r-1}\\
\downarrow & \  & \downarrow\\
Spec\ \overline{S} & \ra & Spec\ \overline{S^r}\end{array}$$
where the bottom horizontal arrow is defined by the surjection
$$\overline{S^r}\ra \frac{\overline{S^r}}{(\overline{f^{\flat}},\overline{\d f^{\flat}},...,\overline{\d^{r-1}f^{\flat}})}=\overline{S},$$
cf. Lemma \ref{becar}.
The natural $\bZ_p^{\times}$-equivariant homomorphism
$M^{\infty}\ra \widehat{\bW_g}$
maps
$f^1,\d f^1,\d^2 f^1,...$
 into $0$; cf. Proposition \ref{funu}. So this homomorphism also maps
$f^{\flat}, \d f^{\flat}, \d^2 f^{\flat},...$
into $0$. One the other hand  this homomorphism also maps
$f^{\partial}-1, \d(f^{\partial}-1), \d^2(f^{\partial}-1),...$
into $0$. So we get an induced $\bZ_p^{\times}$-equivariant homomorphism
$\overline{N^r}\ra \overline{\bW}_g$, hence (by restriction) we get a
$\bZ_p^{\times}$-equivariant homomorphism $\Sigmat^{r,r-1}\ra \overline{\bW}_g$, and hence we get an induced
$\bZ_p^{\times}$-equivariant homomorphism
$$T^r=
\frac{\Sigmat^{r,r-1}}{(\overline{f^{\flat}},...,\overline{\d^{r-1}f^{\flat}})}\ra \overline{\bW}_g.$$
Since $Spec\ \overline{\bW}_g$ is  irreducible the closure $Z$ of the image of $Spec\ \overline{\bW}_g \ra Spec\ T^r$ is contained in one of the connected components
of $Spec\ T^r$. Since $Z$ dominates $Spec\ \overline{S}$ and since $Spec\ T^r$ is finite and \'{e}tale over $Spec\ \overline{S}$, it follows that $Z$ is a connected
component of $Spec\ T^r$.
Note that $Z$ is a $\bZ_p^{\times}$-invariant subset of $Spec\ T^r$, hence $\Gamma_{r-1}$-invariant.
Recall that by (\ref{zty}) $\Gamma_{r-1}$ acts transitively on the fibers of $Spec\ \Sigmat^{r,r-1}\ra Spec\ \overline{S^r}$.
Hence $\Gamma_{r-1}$ acts transitively on the fibers of $Spec\ T^r \ra Spec\ \overline{S}$. Since each connected component of $Spec\ T^r$
surjects onto $Spec\ \overline{S}$ and since $Z$ is $\Gamma_{r-1}$-invariant
it follows  that $Spec\ T^r$ must be connected.
Since $Spec\ T^r$ is  smooth over $k$ and connected it follows that $T^r$ is an integral domain.
Since  $Spec\ T^r$ is connected it must coincide with $Z$ hence $Spec\ \overline{\bW}_g \ra Spec\ T^r$
is dominant. Since  $T^r$ is an integral domain,  $T^r \ra \overline{\bW}_g$ is injective.
So $\lim_{\ra} T^r \ra \overline{\bW}_g$ is injective. But
$$\lim_{\ra} T^r=\lim_{\ra} \Sigmat^{r,r-1}/(\overline{f^{\flat}},\overline{\d f^{\flat}}, \overline{\d^2 f^{\flat}},...)=
\overline{\tS^{\infty}}/(\overline{f^{\flat}},\overline{\d f^{\flat}}, \overline{\d^2 f^{\flat}},...).$$
This proves assertion 3.

On the other hand since each connected component of $Spec\ \Sigmat^{r,r-1}$ surjects onto $Spec\ \overline{S^r}$ and
$Spec\ T^r$ is connected
it follows that
$Spec\ \Sigmat^{r,r-1}$ itself is connected. This ends the proof of the Theorem.
\qed

\begin{corollary}
\label{tzutzu}
Assume $X=Spec\ S$ is a modular Fourier-framed curve  with $E_{p-1}$ invertible on $X$.
Let $f(q)\in R((q))$ be contained in the image of the map
$M\otimes_R K \ra  K((q))$. Then $f(q)$
is contained in the image of the map $M^{\infty}\ra R((q))\hat{\ }$.
\end{corollary}

{\it Proof}.
Write $f(q)=E(\frac{G}{p^{\nu}})=\frac{G(q)}{p^{\nu}}$, where $G \in M$.
The image of $G$ in $R((q))\otimes \bZ/p^{\nu}\bZ$ is $0$ so the image of $G$ in $\tS^{\infty} \otimes \bZ/p^{\nu}\bZ$ is in the kernel of
$\tS^{\infty}\otimes \bZ/p^{\nu}\bZ\ra S^{\infty}_{for}\otimes\bZ/p^{\nu}\bZ$. But the latter morphism is injective; indeed this
is trivially checked by induction on $\nu$, using
 Theorem \ref{babanovac}.
It follows that the image of $G$ in $\tS^{\infty}\otimes \bZ/p^{\nu}\bZ$ is $0$, hence the image of $G$ in $\tS^{\infty}$ belongs to $p^{\nu}\tS^{\infty}$. Hence the image of $G$ in $R((q))\h$ belongs to $p^{\nu}\cdot Im(M^{\infty}\ra R((q))\h)$.
It follows that $f(q)$ belongs to the image of $M^{\infty}\ra R((q))\h$.
\qed

\medskip

Recall that we denoted by $Y_{ord}$ the locus in $Y:=Y_1(N)_R$ where $E_{p-1}$ is invertible.

\begin{corollary}
\label{umbirac}
Consider the modular Fourier-framed curve $X=Spec\ S=Y_{ord}$.  Then
  the image of $M^{\infty}\ra R((q))\hat{\ }$ contains $\bD$ and hence is $p$-adically dense in $\bW$.
\end{corollary}

{\it Proof}.
By Corollary \ref{tzutzu} the image of $M^{\infty}\ra R((q))\hat{\ }$ contains the  ring $\bD$. But this image also contains
the ring $R[\Delta^{-1}]$. We conclude by Lemma \ref{denn}.
\qed

\medskip

Let us also note the following refinement of our results on conjugate derivations.

\begin{proposition}
Assume $X=Spec\ S$ is a modular Fourier-framed curve
and assume the Eisenstein forms $E_4,E_6,E_{p-1}$
are invertible on $X$. Then the sequence
$$\partial^*_0,\partial^*_1,\partial^*_2,...\in Der_R(\tS^{\infty},\tS^{\infty})$$
is a pro-basis of $Der_R(\tS^{\infty},\tS^{\infty})$.
\end{proposition}

{\it Proof}.
Assume $\sum_{j\geq 0} c_j \partial^*_j=0$ with $c_j \in \tS^{\infty}$ and let us prove that $c_j=0$ for all $j$. Assume this is not the case.
We may assume not all $\overline{c}_j$ are $0$ in $\overline{\tS^{\infty}}$.
We get $\sum_j \overline{c}_j \overline{\partial^*_j}=0$ in $Der_k(\overline{\tS^{\infty}},\overline{\tS^{\infty}})$.
By Proposition \ref{noua} and Theorem \ref{babanovac} we get $\overline{c}_j=0$ for all $j$, a contradiction.

Now let $\partial \in Der_R(\tS^{\infty},\tS^{\infty})$. Using Proposition \ref{noua} and Theorem \ref{babanovac} we can find
elements $\gamma_j \in (\tS^{\infty})\hat{\ }$, $j \geq 0$, such that
$$\partial=\sum_{j \geq 0} \gamma_j \partial^*_j:\tS^{\infty} \ra (\tS^{\infty})\hat{\ }.$$
Evaluating $\partial j$ we get $\gamma_0=1$. Since the image of $\partial \d j$ in $(\tS^{\infty})\hat{\ }$ is in
$\tS^{\infty}$ and
$\partial_1 \d j=\phi(\partial j) \in (\tS^{\infty})^{\times}$
we get $\gamma_1\in \tS^{\infty}$. Since the image of $\partial \d^2 j$  in $(\tS^{\infty})\hat{\ }$ is in
$\tS^{\infty}$ and
$\partial_2 \d^2 j=\phi^2(\partial j) \in (\tS^{\infty})^{\times}$
we get $\gamma_2\in \tS^{\infty}$. Continuing in this way we get that $\gamma_j\in \tS^{\infty}$ for all $j$ which ends the proof of the Proposition.
\qed

\section{The  theory for several primes}

In this section several primes will be involved. So it will be important to keep track of the prime $p$
used in the one prime theory of the previous section by making $p$ appear as an index for the various objects considered there. In particular the objects
$$\d, \phi, R,  J^r(X), S^r, \tS^r, M^r, \tM^r, S^r_{for}, f^0, f^1, f^{\partial}, f^{\sharp}, f^{\natural}, H, \Psi, etc.$$
introduced in the previous section in the case of the prime $p$ will be denoted from now on by
 $$\d_p, \phi_p, R_p,  J^r_p(X), S^r_p, \tSp^r, M^r_p, \tMp^r,
S^r_{for,p}, f^0_p, f^1_p, f_p^{\partial}, f^{\sharp}_p, f^{\natural}_p, H_p, \Psi_p, etc.$$

\subsection{Review of concepts and terminology from \cite{laplace}}

For any two distinct rational primes $p_1,p_2$ consider the
polynomial $C_{p_k,p_l}$ in the ring $\bZ[X_0,X_1,X_2]$ defined by
\begin{equation}
\label{commutator}
C_{p_1,p_2}(X_0,X_1,X_2) := \frac{C_{p_2}(X_0^{p_1},p_1X_1)}{p_1}
-\frac{C_{p_1}(X_0^{p_2},p_2 X_2)}{p_2}
 -\frac{\delta_{p_1} p_2}{p_2} X_2^{p_1}+ \frac{\delta_{p_2}p_1}{p_1}
X_1^{p_2}\, .
\end{equation}
Let $\cP=\{p_1,\ldots ,p_d\}$ be a finite set of  primes in $\bZ$.
A $\d_{\cP}$-{\it ring} is a ring $A$  equipped with
$p_k$-derivations $\delta_{p_k}:A \ra A$,  $k=1,\ldots ,d$, such that
\begin{equation}
\label{identit}
\delta_{p_k}\delta_{p_l}a-\delta_{p_l}\delta_{p_k}a=C_{p_k,p_l}(a,
\delta_{p_k}a,\delta_{p_l}a)
\end{equation}
for all $a \in A$, $k,l=1,\ldots ,d$.
A {\it homomorphism
of $\delta_{\cP}$-rings} $A$ and $B$ is a homomorphism
of rings $\varphi: A \rightarrow B$ that commutes
with the $p_k$-derivations in $A$ and $B$, respectively.
If $\phi_{p_k}(x)=x^{p_k}+p_k\d_{p_k}x$ is the homomorphism  associated to $\delta_{p_k}$,
condition
(\ref{identit}) implies that
\begin{equation}
\label{commu}
\phi_{p_k}\phi_{p_l}(a)=\phi_{p_l}\phi_{p_k}(a)
\end{equation}
for all $a \in A$. Conversely, if the commutation relations (\ref{commu})
hold, and the $p_k$s are non-zero divisors in $A$, then conditions
(\ref{identit}) hold, and we have that
$
\phi_{p_k}\delta_{p_l} a  =  \delta_{p_l}\phi_{p_k}a
$
for all $a \in A$. If $A$ is a $\d_{\cP}$-ring then for all $k$, the $p_k$-adic completions
$A^{\widehat{p_k}}$ are $\d_{\cP}$-rings in a natural way.

Note that, since $$Aut(R_{p_k})=Aut(R_{p_k}/p_kR_{p_k})=\lim_{\leftarrow} \bZ/n\bZ$$ is commutative, it follows that to give a $\d_{\cP}$-ring structure on $R_{p_k}$ is equivalent to giving automorphisms $$\phi_{p_1},...,\phi_{p_{k-1}},\phi_{p_{k+1}},...,\phi_{p_d}$$
 of $R_{p_k}$; the automorphism $\phi_{p_k}$
is, of course, uniquely determined. From now on we shall fix, for each $k$ a $\d_{\cP}$-ring structure on $R_{p_k}$.

For a relation between these concepts and the theory of lambda rings we refer to
\cite{borger} and the references therein.

We let $\bZg=\{0,1,2,3, \ldots \}$, and
let $\bZg^d$ be given the product order. We
let $e_k$ be the element of $\bZg^d$ all of whose components are
zero except the $k$-th, which is $1$. We set $e=\sum e_k=(1,...,1)$.
For $i=(i_1,...,i_d) \in \bZg^d$ we set $\cP^i=p_1^{i_1}...p_d^{i_d}$, $\d_{\cP}^i=\d_{p_1}^{i_1}...\d_{\cP}^{i_d}$,
$\phi_{\cP}^i=\phi_{\cP^i}=\phi_{p_1}^{i_1}...\phi_{p_d}^{i_d}$.
We define $W_{\cP}=\bZ[\phi_{p_1},...,\phi_{p_d}]$, the ring of polynomials in the commutative variables $\phi_{p_1},...,\phi_{p_d}$; this ring has a natural order defined by the condition that a polynomial is $\geq 0$ if and only if its coefficients are $\geq 0$. We set $W_{\cP}(r)=\sum_{i \leq r}\bZ \phi_{\cP}^i$. Define $deg:W_{\cP} \ra \bZ$ by $deg(\sum c_n \phi_n)=\sum c_n$.

A {\it $\delta_{\cP}$-prolongation system} $A^*=(A^r)$ is an inductive
system of rings $A^r$ indexed by $r \in \bZg^d$, provided with
transition maps $\varphi_{rr'}:A^r \ra A^{r'}$ for any pair of indices
$r$, $r^{'}$ such that $r \leq r'$, and equipped with $p_k$-derivations
$$
\delta_{p_k}:A^r\ra A^{r+e_k}\, ,
$$
$k=1,\ldots ,d$, such that (\ref{identit}) holds for all $k$, $l$, and
such that
$$
\varphi_{r+e_k,r'+e_k}\circ \delta_{p_k}=\delta_{p_k} \circ
\varphi_{rr'} :A^r \ra A^{r'+e_k}
$$
for  all $r\leq r'$ and all $k$.
A morphism of prolongation systems $A^* \ra B^*$ is a system of ring
homomorphisms
$u^r:A^r \ra B^r$ that commute with the
$\varphi$s and the $\delta$s of $A^*$ and $B^*$, respectively.

Any $\delta_{\cP}$-ring $A$ induces a $\delta_{\cP}$-prolongation system
$A^*$ where $A^r=A$ for all $r$ and $\varphi=$identity. If $A$ is a
$\delta_{\cP}$-ring and $A^*$ is the associated
$\delta_{\cP}$-prolongation system, we say that
a $\delta_{\cP}$-prolongation system $B^*$ is a
{\it $\delta_{\cP}$-prolongation system over $A$}
if it is equipped with a morphism $A^* \ra B^*$.
We have a natural concept of morphism of $\delta_{\cP}$-prolongation
systems over $A$.

Consider the ring  $$\bZ_{(\cP)}=\bigcap_{k=1}^d \bZ_{(p_k)} \subset \bQ.$$
For any affine scheme of finite type $X$ over $\bZ_{(\cP)}$ we considered in \cite{laplace}
a system of schemes of finite type, $\cJ^r_{\cP}(X)$ over $\bZ_{(\cP)}$, called the $\d_{\cP}$-jet spaces of $X$; if $X=Spec\ \bZ_{(\cP)}[x]/(f)$ then
$$\cJ^r_{\cP}(X)=Spec\ \bZ_{(\cP)}[\d_{\cP}^ix;i\leq r]/(\d_{\cP}^if;i\leq r).$$
Cf. also \cite{borger}, where these spaces were introduced independently.
The rings $\cO(\cJ^r_{\cP}(X))$ form then a $\d_{\cP}$-prolongation system.

 \subsection{Complements to \cite{laplace}}

 \subsubsection{Splitting of completions}
 We will need the following {\it splitting result}
 for the various $p_k$-adic completions of the $\d_{\cP}$-jet spaces; here $X$ is an affine scheme of finite type over $\bZ_{(\cP)}$ .

\begin{proposition}
We have a natural isomorphism of $\d_{\cP}$-prolongation systems
$$\cO(\cJ_{\cP}^r(X))^{\widehat{p_k}}\simeq \left(
\bigotimes_{i \leq r-r_ke_k} \cO(\cJ_{\{p_k\}}^{r_k}(X_i)^{\widehat{p_k}}) \right)^{\widehat{p_k}},$$
where, in the right hand side,  $\otimes=\otimes_{\bZ_{p_k}}$, $X_i=X$ for all $i$, $\d_{p_k}$ is induced by the operators
$$\d_{p_k}:\cO(\cJ_{\{p_k\}}^{r_k}(X_i))\ra \cO(\cJ_{\{p_k\}}^{r_k+1}(X_i)),$$
and $\d_{p_j}$ ($j \neq k$) are defined by the automorphisms $\phi_{p_j}$ induced from
 $$\cO(\cJ_{\{p_k\}}^{r_k}(X_i))\simeq \cO(\cJ_{\{p_k\}}^{r_k}(X_{i+e_j})).$$
\end{proposition}

{\it Proof}. One shows that the right hand side satisfies the universality property of the left hand side;
the universality property in question is explained in \cite{laplace}, Remark 2.20.

\begin{corollary}
\label{mordeso}
We have a natural isomorphism of $\d_{\cP}$-prolongation systems
$$\left(\cO(\cJ_{\cP}^r(X))\otimes_{\bZ_{(\cP)}}  R_{p_k}\right)^{\widehat{p_k}}\simeq \left(
\bigotimes_{i \leq r-r_ke_k} \cO(J_{p_k}^{r_k}(X_{i,k})) \right)^{\widehat{p_k}},$$
where, in the right hand side, $\otimes=\otimes_{R_{p_k}}$ and $X_{i,k}=X_{R_{p_k}}$ for all $i$.
\end{corollary}

\subsubsection{Analytic continuation}
\label{anacon}
Let us explain, in an abstract setting,  a concept that will play a key role
in the concrete discussion of the next section. This concept is a generalization
of a concept introduced in \cite{laplace}.

Assume we are given a set of primes $\cP=\{p_1,...,p_d\}$ and consider the $\d_{\cP}$-prolongation system
\begin{equation}
\label{sfork}
\begin{array}{rcl}
\cS^r_{for,k} & := & R_{p_k}[[\d^i_{\cP}t;i \leq r, i_k=0]][\d^i_{\cP}t;i \leq r, i_k\geq 1]^{\widehat{p_k}}\\
\  & \  & \  \\
\  & \simeq & R_{p_k}[[t_i;\ i \leq r-r_ke_k]][\d_{p_k}^j t_{i};\ i \leq r-r_ke_k,1 \leq j \leq r_k]^{\widehat{p_k}},
\end{array}
\end{equation}
where $t_i$ correspond to $\phi_{\cP}^it$ for $i \leq r-r_ke_k$.

Assume now that:

1) For each $k=1,...,d$ we are given a $\d_{\cP}$-prolongation system
\begin{equation}
\label{As}
(\cA^r_k);
\end{equation}

2) For each $k$ we are given a morphism of $\d_{\cP}$-prolongation systems
\begin{equation}
\label{Es}
\cE_{k}^r:\cA^r_k \ra \cS^r_{for,k}.
\end{equation}

Finally consider the ring
\begin{equation}
\cS^r_{for,0}:=\bZ_{(\cP)}[[\d_{\cP}^it;i \leq r]].\end{equation}
For each $k$ we shall view both $\cS^r_{for,k}$ and $\cS^r_{for,0}$  as  subrings of
$R_{p_k}[[\d_{\cP}^i t;i \leq r]]$.

\begin{definition}
We say that a  family
$$(f_1,...,f_d) \in  \prod_{k=1}^d \cA^r_k$$
 can be {\it analytically continued} if  there exists (a necessarily unique)
$f_0 \in \cS^r_{for,0}$ such that for all $k=1,...,d$ we have
$$\cE^r_k(f_k)=f_0\in  R_{p_k}[[\d_{\cP}^i t;i \leq r]].$$
We say that $(f_1,...,f_d)$ has {\it expansion} $f_0$.
\end{definition}

\subsection{Application to $\d$-modular forms}

In what follows we specialize our discussion to the case of differential modular forms.

\subsubsection{The rings $\cS^r_k,\cM^r_k$}
We assume $Y_{\bZ[1/N]}$ and  $Y_{\bZ[1/m]}$ are as in sections \ref{modu} and \ref{shimu} respectively.
In the first case we assume $N$ is not divisible by any of the primes in $\cP$,
and in the second case we assume the primes in $\cP$ are  sufficiently big.
We let $Y_{\bZ_{(\cP)}}$ and $Y_{R_{p_k}}$ be the curves over $\bZ_{(\cP)}$ and $R_{p_k}$ respectively obtained via base change.
As in sections \ref{modu} and \ref{shimu} we let $X_{R_{p_k}} \subset Y_{R_{p_k}}$ be open affine subsets.
We assume that $\overline{X_{R_{p_k}}}$ is contained in the ordinary locus. We also assume, for simplicity,
that $\overline{X_{R_{p_k}}}$ is principal in $\overline{Y_{R_{p_k}}}$ defined by a function $s_k \in \cO(\overline{Y_{R_{p_k}}}$)
and that $L_{R_{p_k}}$ is trivial on $X_{R_{p_k}}$ with basis $x_k$.
(Although the choice of $X_{R_{p_k}}$ is arbitrary with the above properties, the theory is easily seen to be independent, in an ``obvious sense" from these choices.)

Now  define
\begin{equation}
\label{mariaa}
\begin{array}{rcl}
\cS^r_k & := & ((\cO(\cJ^r_{\cP}(Y_{\bZ_{(\cP)}}))\otimes_{\bZ_{(\cP)}} R_{p_k})[\phi_{\cP}^is_k^{-1}; i \leq r])^{\widehat{p_k}}\\
\  & \  & \  \\
\  & = & \left( \bigotimes_{i \leq r-r_ke_k}\cO(J_{p_k}^{r_k}(X_{i,k}))\right)^{\widehat{p_k}},\end{array}\end{equation}
where $X_{i,k}=X_{R_{p_k}}$ for all $i$.
(The last equality follows from Corollary \ref{mordeso}.)
Then $(\cS^r_k)$ has a natural structure of $\d_{p_k}$-prolongation sequence. Actually,
 once we have fixed a $\d_{\cP}$-ring structure on $R_{p_k}$, $(\cS^r_k)$ has a naturally induced structure of $\d_{\cP}$-prolongation system (extending the previous one). We also have a natural morphism of $\d_{p_k}$-prolongation sequences $S^{r_k}_{p_k}\ra \cS^r_k$.
Set
\begin{equation}
\label{katie}
\begin{array}{rcl}
\cM^r_k& := & \cS^r_k[\d_{\cP}^ix_k,\phi_{\cP}^i x_k^{-1};i \leq r]^{\widehat{p_k}},\\
\  & \  & \  \\
\  & = & \left( \bigotimes_{i \leq r-r_ke_k}\cO(J_{p_k}^{r_k}(X_{i,k}))[\d_{p_k}^jx_{i,k}, x_{i,k}^{-1};0 \leq j\leq r_k]^{\widehat{p_k}}\right)^{\widehat{p_k}}\\
\  & \  & \  \\
\cM^r_k(w) & := & \cS^r_k x^w_k,
\end{array}
\end{equation}
where $x_{i,k}$ are new variables, ``copies" of $x_k$.
For each $k$, $(\cM^r_k)$ is a $\d_{\cP}$-prolongation system and we have a natural morphism of $\d_{p_k}$-prolongation sequences $M^{r_k}_{p_k}\ra \cM^r_k$. In particular we may consider the image
of $f^1_{p_k}\in M^1_{p_k}(-1-\phi_{p_k})$ in $\cM^{e_k}_k(-1-\phi_{p_k})$ which we continue to denote by $f^1_{p_k}$. Similarly, if
 $f =\sum a_n q^n$
is a newform of weight $2$ on $\Gamma_0(N)$ with $\bZ$-Fourier coefficients and if the primes in $\cP$ are sufficiently big
then we may consider the image of $f^{\sharp}_{p_k}\in M^2_{p_k}(0)$ in $\cM^{2e_k}_k(0)$ which we continue to denote by $f^{\sharp}_{p_k}$.

Furthermore define
$$\begin{array}{lll}
\cS^{\infty}_k & := & \lim_{\ra}\cS^r_k,\\
\  & \  & \  \\
\cM^{\infty}_k&:=&\lim_{\ra}\cS^r_k,\\
\  & \  & \  \\
 \cS^{\infty}_{for,k}&:=&\lim_{\ra}\cS^r_{for,k}.\end{array}$$

 \subsubsection{$\d_{\cP}$-expansion maps}
 Next we will construct  a natural
morphism of $\d_{\cP}$-pro\-long\-ation systems $$\cE^r_{k}:\cM^{r}_k \ra \cS^r_{for,k},$$
called $\d$-{\it expansion maps},   as follows.

Assume first we are in the ``Fourier case".
Consideration of the Tate curve yields a homomorphism
$\cO(Y_{\bZ_{(\cP)}}) \ra R_{p_k}((q))^{\widehat{p_k}}$. Composing
this with the homomorphism
$$R_{p_k}((q))^{\widehat{p_k}}\ra R_{p_k}[[t]],\ \ \ q\mapsto t+1$$
we get a homomorphism
$\cO(Y_{\bZ_{(\cP)}}) \ra R_{p_k}[[t]]$.
By universality we obtain a homomorphism of $\d_{\cP}$-prolongation systems
$$\cO(\cJ^r_{\cP}(Y_{\bZ_{(\cP)}}))\ra \cS^r_{for,k}.$$
Using the fact that $E(s_k)$ is invertible in $R_{p_k}[[t]]$
we get a homomorphism
$\cS^r_k \ra \cS^r_{for,k}$. Sending $x_k \mapsto \cE(x_k)$
 we get the desired homomorphism $\cE^r_k:\cM^{r}_k \ra \cS^r_{for,k}$.
Note  that this homomorphism extends, in the obvious sense, the homomorphism
$$M^{r}_{p_k}\stackrel{\cE^r}{\longrightarrow} S^r_{for,p_k}\stackrel{\iota}{\longrightarrow} R_{p_k}[[t]][\d_{p_k} t,...,\d_{p_k}^rt]^{\widehat{p_k}},$$
where $\iota$ is induced by $q \mapsto t+1$. Consequently, by Proposition \ref{funu}, we have
\begin{equation}
\label{expn2}
\cE^{e_k}_k(f_{p_k}^1)=\Psi_{p_k}.
\end{equation}

Next we consider the ``Serre-Tate" case. For this case we assume we are given a collection of points
\begin{equation}
\label{Pk}
(P_k)_{1\leq k\leq d},
\end{equation}
$P_k \in X_{R_{p_k}}(R_{p_k})$, as in sections \ref{modu} and \ref{shimu} respectively. (We do not assume any compatibility between the points $P_k$. One of the most remarkable aspects of the theory is that
its constructions and results are valid without assuming any compatibility between the $P_k$s. Also the theory is, in a sense easily made precise, independent of the choice of the $P_k$s.) Using the decompositions
(\ref{sfork}) and (\ref{mariaa}) plus the one prime theory we get induced maps $E^r_k:M^r_k \ra S_{for,k}^r$. Once again (\ref{expn2}) holds.

Note  the following $\d_{\cP}$-{\it expansion principle}:

\begin{proposition}
\label{pissed}
The homomorphism $$\overline{\cE^r_k}:\overline{\cS^r_k}\ra \overline{\cS^r_{for,k}}$$ is injective.
In particular the homomorphisms
$$\cE^r_{k}:\cM^{r}_k(w)\ra \cS^r_{for,k}$$
are injective with torsion free cokernel.
\end{proposition}

{\it Proof}.
This follows from the one prime situation (cf. Proposition \ref{expansionprinciple}) plus the ``splittings"
(\ref{sfork}) and (\ref{mariaa}).
\qed

\subsubsection{The rings $\tSk^r,\tMk^r$}

As in the case of one prime we define
$$\begin{array}{rcl}
\tSk^r & := & Im(\cE^r_k:\cM^r_k\ra \cS^r_{for,k})\\
\  & \  & \  \\
\tMk^r & := & (\tSk^r \otimes_{\cS^r_k} \cM^r_k)^{\widehat{p_k}},\\
\  & \  & \  \\
\tSk^{\infty} & = & \lim_{\ra} \tSk^r,\\
\  & \  & \  \\
  \tMk^{\infty} & = & \lim_{\ra} \tMk^r.\end{array}$$
The inclusion $\tSk^r \subset \cS^r_{for,k}$ and the homomorphism $\cE^r_k:\cM^r_k \ra \cS^r_{for,k}$ induce a homomorphism (still denoted by $\cE^r_k$ and still referred to as {\it $\d_{\cP}$-expansion map}),
$$\cE^r_k:\tMk^r \ra \cS^r_{for,k}.$$
Note that  $\tcS^r$ is $p_k$-adically complete and if $x_k$ is a basis of $L=L_{R_{p_k}}$ then
$$\tMk^r=\tSk^r[\d_{\cP}^ix_k,\phi_{\cP}^ix_k^{-1};i \leq r]^{\widehat{p_k}}.$$
Note that $(\tSk^r)$, $(\tMk^r)$ have natural structures of $\d_{\cP}$-prolongation systems.
Define
$$\tMk^r(w):=\tSk^r x_k^w \subset \tMk^r;$$
the latter definition is independent of the choice of the basis $x_k$.
We have  morphisms of $\d_{p_k}$-prolongation sequences $\tSpk^r\ra \tSk^r$, $\tMpk^{r_k}\ra \tMk^r$.
In particular we may consider the image of $f^0_{p_k}\in \tMpk^0(1)$ into $\tMk^0(1)$ which we will still denote by $f^0_{p_k}$. Clearly we have $\cE^0_k(f^0_{p_k})=1$.

\begin{proposition}
\label{lotsofinjcal}
\

1) The homomorphisms
$\overline{\cS^r_k}\ra \overline{\tSk^r}$, $\overline{\cS^{\infty}_k}\ra \overline{\tSk^{\infty}}$, $\overline{M_k^r}\ra \overline{\tMk^r}$, $\overline{\cM^{\infty}_k}
\ra \overline{\tMk^{\infty}}$ are injective. In particular the homomorphisms $\cS^r_k\ra \tSk^r$, $\cS^{\infty}_k\ra \tSk^{\infty}$, $\cM^r_k\ra \tMk^r$, $\cM^{\infty}_k\ra \tMk^{\infty}$
 are injective with torsion free cokernel.

2) The homomorphisms $\tMk^r(w)\ra \cS^r_{for,k}$ are injective.

3) The homomorphisms $\tSk^r\ra \tSk^{r+1}$, $\tMk^r \ra \tMk^{r+1}$ are injective.
\end{proposition}

{\it Proof}. Use Proposition \ref{lotsofinj} and the ``splittings" (\ref{sfork}) and (\ref{mariaa}).
\qed

\begin{remark}
 Again $\overline{\tSk^{\infty}}$ is not a priori an integral domain; but it possesses a natural quotient which is an integral domain:
$$\stuffk:=Im(\overline{M^{\infty}_k}\ra \overline{S^{\infty}_{for,k}}).$$
\end{remark}

By the proof of Theorem \ref{main} and the ``splittings" (\ref{sfork}) and (\ref{katie}) we get:

\begin{theorem}
\label{maincal}
The ring $\overline{\tSk^{\infty}}$ is a quotient of an ind-\'{e}tale   $\Gamma$-extension of  $\overline{\cS^{\infty}_k}$, where
$\Gamma$ is a profinite abelian group.
\end{theorem}

\begin{corollary}\

1) $\overline{\tSk^{\infty}}$ is an integral extension of  $\overline{\cS^{\infty}_k}$.

2) $\stuffk$ is an ind-\'{e}tale $\Gamma'$-extension of $\overline{\cS^{\infty}_k}$, where $\Gamma'$ is a closed subgroup of $\Gamma$.
\end{corollary}

 Theorem \ref{maincal}
can be morally viewed as saying that the  ``$\d_{\cP}$-Igusa curve"
(of which the rings $\tSk^{\infty}$ are an incarnation) is  a ``formal profinite cover"
(embedding into an ``abelian formal pro-\'{e}tale cover")
 of the modular/Shimura curve (whose $\d_{\cP}$-geometric incarnation are the rings $\cS^{\infty}_k$).

One can prove analogues, in our several prime setting here, of Theorem \ref{babanovac} and Propositions \ref{zece}, \ref{noua}.
We are not going to need these analogues so we are not going to state them explicitly. Instead, we will concentrate, in what follows,
 on the purely ``several prime concept''
of analytic continuation.

In the following definition the rings $M^r_k$ will play the role of the rings $A^r_k$ in (\ref{As})
and the $\d_{\cP}$-expansion maps $E^r_k:M^r_k \ra S^r_{for,k}$ will play the role of the maps (\ref{Es}).

\begin{definition}
\label{igsttt}
A $\d_{\cP}$-{\it modular form of weight $w \in W_{\cP}$ and order $r \in \bZg^d$} is a family
$$f=(f_1,...,f_d) \in \prod_{k=1}^d M_k^r$$
 that can be analytically continued and such  that $f_k \in M^r_k(w)$ for all $k$. We denote by $M^r_{\cP}(w)$
the group of all such forms. There is a naturally induced {\it expansion map}
$$\cE^r:M^r_{\cP}(w) \ra \cS^r_{for,0}$$
which by Proposition \ref{lotsofinjcal} is injective.
\end{definition}

Similarly we may take the rings $\tMk^r$ to play the role of the rings $\cA^r_k$ in (\ref{As})
and the $\d_{\cP}$-expansion maps $\cE^r_k:\tMk^r \ra \cS^r_{for,k}$ to play the role of the maps (\ref{Es}).

\begin{definition}
\label{igst}
An {\it Igusa $\d_{\cP}$-modular form of weight $w \in W_{\cP}$ and order $r \in \bZg^d$} is a family
$$f=(f_1,...,f_d) \in \prod_{k=1}^d \tMk^r$$
 that can be analytically continued and such  that $f_k \in \tMk^r(w)$ for all $k$. We denote by $\tMP^r(w)$
the group of all such forms. There is a naturally induced {\it expansion map}
$$\cE^r:\tMP^r(w) \ra \cS^r_{for,0}$$
which by Proposition \ref{lotsofinjcal} is injective.
\end{definition}

So $M_{\cP}^r(w)$ is a $\bZ_{(\cP)}$-submodule of $\tMP^r(w)$.

\begin{remark} If, in the above definitions, we are in the Fourier expansion case then we view $f$ as being ``analytically continued along the section $\infty$".
If we are in the Serre-Tate expansion case and if the collection of points $(P_k)$ in (\ref{Pk}) comes from a $\overline{\bQ}$-point $P$ of the modular/Shimura curve, then we view $f$ as being ``analytically continued along $P$". Here, we say that $(P_k)$ {\it comes from}
 $P$ if
 for each $k$ the elliptic (respectively false elliptic) curve corresponding to $P_k$ is isomorphic over $\overline{\bQ}$ to the  curve corresponding to $P$.
\end{remark}

\subsubsection{Basic examples: the forms $f^0$, $f^e$, $f^{2e}$}
We may now introduce some of the fundamental objects of our ``several primes" theory.

First recall that we have at our disposal forms $f^0_{p_k}\in \tMk^0(1)$. Then we have the following obvious

\begin{proposition}
The family
$$f^0:=(f^0_{p_1},...,f^0_{p_d})\in \prod_{k=1}^d \tMk^0$$
is an Igusa $\d_{\cP}$-modular form of weight $1$ and order $0$; i.e.
$f^0\in \tMP^0(1)$. Moreover $f^0$   has expansion  $\cE^0(f^0)=1$.
\end{proposition}

In particular, for any $w \in W_{\cP}(r)$ we have
$$(f^0)^w=((f^0_{p_1})^w,...,(f^0_{p_d})^w) \in \tMP^r(w).$$
Note that if the weight $w$ is divisible in $W_{\cP}$ by $(\phi_{p_1}-1)...(\phi_{p_d}-1)$
then  $(f^0)^w$ is actually belongs to $\tMP^r(w)$.

Next we introduce a form of order
$e:=(1,...,1)$
which we shall call $f^e$. Indeed recall the forms $f^{\natural}_{p_k}$ and set, for each $k=1,...,d$,
\begin{equation}
\label{fimportant}
f^e_{k}  :=(-1)^{d-1} \left( \prod_{l \in I_k} \left(1-\frac{\phi_{p_l}}{p_l}\right)
 \right) f^{\natural}_{p_k} \in \tSk^r=\tMk^r(0),
\end{equation}
where $I_k=\{1,...,d\}\backslash \{k\}$.
Also set
\begin{equation}
\label{zerolevel}
f^e_0:=\frac{1}{p_1...p_d}(\phi_{p_1}-p_1)...(\phi_{p_d}-p_d)\log (1+t) \in \cS^e_{for,0}=\bZ_{(\cP)}[[\d_{\cP}^i t;i \leq e]].
\end{equation}
The fact that the above series belongs to $\bZ_{(\cP)}[[\d_{\cP}^i t;i \leq e]]$ (and not merely to
$\bQ[[\d_{\cP}^i t;i \leq e]]$) is due to the fact that
$$\frac{1}{p_k}(\phi_{p_k}-p_k)\log (1+t)\in \bZ_{p_k}[[t]][\d_{p_k}t]^{\widehat{p_k}}$$
for each $k$.

\begin{theorem}
\label{docct}
The family
$$f^e:=(f^e_{1},...,f^e_{d}) \in \prod_{k=1}^d \tMk^e$$
is an Igusa $\d_{\cP}$-modular form of weight $0$ and order $e$; i.e.
$f^e \in \tMP^e(0)$. Moreover $f^e$   has expansion $f^e_0$; i.e. $\cE^e(f^e)=f^e_0$.
\end{theorem}

{\it Proof}.
Recall that
$\cE^0_k(f^0_{p_k})=1$.
On the other hand, by (\ref{expn2}),
$$\cE^e_k((f^1_{p_k})^{\phi_{n}}) =\Psi_{p_k}^{\phi_{n}}=\frac{1}{p_k}\phi_{n}(\phi_{p_k}-p_k)\log (1+t).$$
Consequently
$$
\cE^e_k(f^e_{k})  = (-1)^{d-1}
 \left( \prod_{l \in I_k} \left(1-\frac{\phi_{p_l}}{p_l}\right)
 \right) \frac{1}{p_k}(\phi_{p_k}-p_k)\log (1+t)= f^r_0.$$
\qed

\begin{remark}
Using the forms $f^0$ and $f^e$ we may construct, for any $w \in W_{\cP}(r)$, with $r \geq e$, the Igusa $\d_{\cP}$-modular form
$$(f^0)^w(f^e)^{\phi_{\cP}^{r-e}} \in \tMP^r(w).$$
\qed
\end{remark}

Finally we introduce  $\d_{\cP}$-modular forms of order $2e$ and weight $0$  which will be called $f^{2e}$. Indeed assume $f =\sum a_n q^n$
is a newform of weight $2$ on $\Gamma_0(N)$ with $\bZ$-Fourier coefficients. Assume, for simplicity, that $f$ is not of CM type. (There is an analogue of what follows for the CM type case.)
 Assume the primes in $\cP$ are sufficiently big.
Define
$$
f^{2e}_k = \left( \prod_{l \in I_k} \left(
 1-a_{p_l}\frac{\phi_{p_l}}{p_l}+p_l\left( \frac{\phi_{p_l}}{p_l}\right)^2
 \right)\right)  f^{\sharp}_{p_k}\in \cM^{2e}_k(0),
$$
where $I_k=\{1,\ldots ,d\}\backslash \{k\}$.
Also set
\begin{equation}
f^{2e}_0=\frac{1}{p_1...p_d} (\phi_{p_1}^2-a_{p_1} \phi_{p_1}+p_1)...
(\phi_{p_d}^2-a_{p_d} \phi_{p_d}+p_d) \sum_n \frac{a_n}{n} q^n \in \cS^{2e}_{for,0}.\end{equation}
Exactly as in the case of Theorem \ref{docct} we get

\begin{theorem}
The family
$$f^{2e}:=(f^{2e}_1,...,f^{2e}_d) \in \prod_{k=1}^d M_k^{2e}$$
is a $\d_{\cP}$-modular form of weight $0$ and order $2e$; i.e.
$f^{2e} \in \cM^{2e}_{\cP}(0)$. Moreover $f^{2e}$   has expansion $f^{2e}_0$; i.e. $\cE^{2e}(f^{2e})=f^{2e}_0$.
\end{theorem}

\subsubsection{Isogeny covariant Igusa $\d_{\cP}$-modular forms}

The following definition extends the concept of isogeny covariance introduced in \cite{difmod, shimura}.
For each $k=1,...,d$ let us fix, once and for all an element  $\gamma_k \in \bZ_{p_k}$
such that $\gamma_k \not\equiv 0,1$ mod $p_k$ and such that $\gamma_k$ is not a root of unity;
the theory that follows is, in a sense that can be made precise, independent of the choice of $\gamma_k$.
Set
$$[\gamma_k](t):=(1+t)^{\gamma_k}-1\in \bZ_{p_k}[[t]].$$

\begin{definition} A series
$$F=F(...,\d_{\cP}^i t,...) \in \cS^r_{for,k},$$
is called {\it isogeny covariant} of degree $\nu \in \bZ$ if
$$F(...,\d_{\cP}^i([\gamma_k](t)),...)=\gamma_k^{\nu} F(...,\d_{\cP}^i t,...).$$
Let $w \in W_{\cP}(r)$ have even degree $deg(w)$.
An element $f_k \in M^r_k(w)$ (respectively an element $f_k\in \tMk^r(w)$)  is called {\it isogeny covariant}
if its expansion $\cE^r_k(f_k)\in \cS^r_{for,k}$ is isogeny covariant of degree $-\frac{deg(w)}{2}$.
An $\d_{\cP}$-modular form $f=(f_k)\in M^r_{\cP}(w)$ (respectively an Igusa  $\d_{\cP}$-modular form  $f=(f_k) \in \tMP^r(w)$) is called {\it isogeny covariant}
if $f_k$ is isogeny covariant for all $k$.
 We denote by $I^r_{\cP}(w)$ (respectively by
$\tIP^r(w)$) the $\bZ_{(\cP)}$-module of all isogeny covariant  $\d_{\cP}$-modular forms in $M^r_{\cP}(w)$ (respectively the $\bZ_{(\cP)}$-module of all isogeny covariant Igusa  $\d_{\cP}$-modular forms in $\tMP^r(w)$).
Hence
$$\cI^r_{\cP}(w):=\tIP^r(w) \cap \cM^r_{\cP}(w)\subset \tMP^r(w).$$
\end{definition}

Note that  the forms $f^0_{p_k}\in \tMk^0(1)$ and $f^1_{p_k}\in \tMk^{e_k}(-1-\phi_{p_k})$ are isogeny covariant. So we have:

\begin{corollary}
Assume $w \in W_{\cP}(r)$, $deg(w)=-2$, $r \geq e$. Then
the form $(f^0_{p_k})^w (f^e_{k})^{\phi_{\cP}^{r-e}} \in \tMk^r(w)$ is isogeny covariant. In other words
the form $(f^0)^w (f^e)^{\phi_{\cP}^{r-e}}$ is isogeny covariant; i.e.  $(f^0)^w (f^e)^{\phi_{\cP}^{r-e}} \in \tIP^r(w)$.
\end{corollary}

Similarly we have:

\begin{corollary}
Assume $w \in W_{\cP}(r)$, $deg(w)=0$. Then the form $(f^0)^w$ is isogeny covariant; i.e.
$(f^0)^w \in \tIP^r(w)$. Moreover, if $w$ is divisible by $(\phi_{p_1}-1)...(\phi_{p_d}-1)$ in $W_{\cP}$ then $(f^0)^w \in I_{\cP}^r(w)$.
\end{corollary}

\begin{remark}
Isogeny covariance is a property which is stronger than the property of being a ``$\d_{\cP}$-Hecke eigenform";
cf. \cite{difmod} for a discussion of this in the case of one prime.
The $\d_{\cP}$-modular form $f^{2e}$ is not isogeny covariant. Nevertheless, by what was shown in \cite{eigen} $f^{2e}$ is, in an appropriate sense, a ``$\d_{\cP}$-Hecke eigenform".
\end{remark}

 \subsubsection{Structure of $\tIP^r(w)$ for $deg(w)=-2,0$}

 Here is the main result of this second part of our paper. It is a structure theorem for the module of isogeny covariant
 Igusa $\d_{\cP}$-modular forms of any order and any weight $w$ with $deg(w)=-2$.

 \begin{theorem}
 \label{uniqueness}
 Let $w \in W_{\cP}(r)$ with $deg(w)=-2$, $r\geq e$. Then
 the $\bZ_{(\cP)}$-module of isogeny covariant $\d_{\cP}$-modular forms $\tIP^r(w)$ is free of rank $r_1...r_d$ with basis
 $$\{(f^0)^w (f^e)^{\phi_{\cP}^{s-e}}\ ;\ e \leq s \leq r\}.$$
 \end{theorem}

 We need a couple of Lemmas.

 \begin{lemma}
\label{funtional}
Assume $F\in \cS^r_{for,k}$ is isogeny covariant of degree $1$ and let $K_{p_k}=R_{p_k}[1/p_k]$.
Then the image of $F$ in $K_{p_k}[[\d_{\cP}^it;i\leq r]]$
 is a $K_{p_k}$-linear combination of the series
$$\phi_{\cP}^i(\log(1+t)),\ \ i\leq r.$$
\end{lemma}

{\it Proof}.
Same argument as in the proof of \cite{shimura}, Lemma 3.9.
\qed

 \begin{lemma}
 \label{ajut}
 Consider the action $\star$ of $\bQ[\phi_{p_1},...,\phi_{p_d}]$ on
 $\bZ_{p_k}[[t]]\otimes \bQ$ defined by
  $\phi_n \star t:=t^n$.
If a polynomial $\Lambda=\sum \lambda_n \phi_n \in
\bQ[\phi_{p_1},... ,\phi_{p_d}]$
satisfies
$$
\Lambda \star \log(1+t) \in \bZ_{p_k}[[T]]\otimes \bQ
$$
for some $k\in \{1,...,s\}$, then $\Lambda$ is divisible in
the ring  $\bQ[\phi_{p_1},... ,\phi_{p_s}]$ by $\phi_{p_k}-p_k$.
\end{lemma}

{\it Proof}. This was proved in the course of the proof of Theorem 3.4 in \cite{laplace}; cf. Claim 2 of that proof.
\qed

 {\it Proof of Theorem \ref{uniqueness}}.
 Let $f=(f_1,...,f_d) \in \tMP^r(w)$ be isogeny covariant and fix an index $k$.
 By Lemma \ref{funtional} the expansion
 $\cE^r_{k}(f_k)$, viewed as an element of $K_{p_k}[[\d_{\cP}^it;i \leq r]]$,
 can be written as
 $$\cE^e_{k}(f_k)=(\sum_{n|\cP^r} c_n \phi_n) \log(1+t),$$
 with $c_n \in K_{p_k}$.
 This series, being equal to the expansion $f_0$ of $f$ also belongs to $\bZ_{(\cP)}[[\d_{\cP}^it;i \leq r]]$. Setting $\d_{\cP}^it=0$ for all $i\neq 0$ we get that
  $$(\sum_{n|\cP^r} c_n \phi_n) \star \log(1+t) \in \bZ_{(\cP)}[[t]].$$
  One checks by induction that $c_n \in \bQ$ for all $n$. Now by Lemma \ref{ajut} it follows
  that
  $$\sum_{n|\cP^r}c_n \phi_n=\frac{1}{p_1...p_d}\left(\sum_{e \leq s \leq r} b_{\cP^s} \phi_{\cP^s}\right)(\phi_{p_1}-p_1)...(\phi_{p_d}-p_d)$$
   with $b_{\cP^s} \in \bQ$.
  Hence, for some sufficiently divisible
  integer $N \in \bZ$, the expansion  of $Nf$ equals the expansion of
  $$g:=N\sum_{e \leq s \leq r}b_{\cP^s} (f^0)^w (f^e)^{\phi_{\cP}^{s-e}}.$$
   By the injectivity
of the expansion map $\cE^r:\tMP^r(w) \ra \cS^r_{for,0}$ (cf. Definition \ref{igst}) it follows that
$Nf= g$. By induction (looking at the coefficient of $t$ to the lowest power)  we get $b_{\cP^s} \in \bZ_{(\cP)}$  for all $s$ and we are done.
 \qed

Similarly (and indeed with a simpler argument which we leave to the reader)
 we get the following structure theorem for the module of isogeny covariant
 Igusa $\d_{\cP}$-modular forms of any order and any weight $w$ with $deg(w)=0$.

\begin{theorem}
Let $w \in W_{\cP}(r)$ with $deg(w)=0$. Then
 the $\bZ_{(\cP)}$-module of isogeny covariant $\d_{\cP}$-modular forms $\tIP^r(w)$ is free of rank one with basis
 $(f^0)^w$.
\end{theorem}

\subsubsection{Vanishing of $\cI^e_{\cP}(-2)$}

``Linear arithmetic partial differential operators" exist, as we have seen, on the ``$\d_{\cP}$-Igusa
curves"; but we do not expect them to exist on the modular curve itself. In other words
we expect that there are no non-zero isogeny covariant $\d_{\cP}$-modular forms of weight $w$ with $deg(w)=-2$; i.e., for such $w$s, $\cI^r_{\cP}(w)=0$.
We can prove this in the modular curve case, in the ``simplest case" $d=2$, $r=e$, $w=-2$:

\begin{theorem}
\label{buniqueness}
Assume we are in the modular curve case and assume  $d=2$. Then $\cI^e_{\cP}(-2)=0$.
\end{theorem}

{\it Proof}.
Assume $\cI^e_{\cP}(-2)$ contains a non-zero element $f$. By
 Theorem \ref{uniqueness} we have $f=c (f^0)^{-2}f^e$ for some $c \in \bZ_{(\cP)}$.
 In particular we have
 $$-c(f^0_{p_1})^{\phi_{p_1}-1}f^1_{p_1}+\frac{c}{p_2}(f^0_{p_1})^{\phi_{p_2}+\phi_{p_1}\phi_{p_2}-2}
 (f^1_{p_1})^{\phi_{p_2}} \in \cM^e_1.$$
The first term in the sum above is in $\cM^e_1$ hence the second term must also be in $\cM^e_1$.
Since by Proposition \ref{lotsofinjcal} $\cM^e_1 \ra \tMunu^e$ has torsion free cokernel
it follows that
$$(f^0_{p_1})^{\phi_{p_2}+\phi_{p_1}\phi_{p_2}-2}
 (f^1_{p_1})^{\phi_{p_2}} \in \cM^e_1$$
 and hence
 $$G:=(f^0_{p_1})^{2\phi_{p_2}-2}
 (f^1_{p_1})^{\phi_{p_2}} \in \cM^e_1.$$
 Reducing modulo $p_1$ and raising to power $\frac{p_1-1}{2}$ we get, using Corollary \ref{partdeux}, that
 $$\overline{H}_{p_1}^{\phi_{p_2}-1} ((\overline{f^1_{p_1}})^{\phi_{p_2}})^{\frac{p_1-1}{2}}=
 \overline{G}^{\frac{p_1-1}{2}}$$
 in $\overline{\tMunu^e}$ and hence in $\overline{\cM^e_1}$.
 Writing $G=g x^{\phi_{p_2}-\phi_{p_1p_2}-2}$ with $g \in \cS_1^e$, $\overline{H}_{p_1}=\overline{h}_{p_1} x^{p_1-1}$ with $\overline{h}_{p_1} \in \overline{\cS^0_1}$, and $f^1_{p_1}=\eta x^{-1-\phi_{p_1}}$ with
 $\eta\in \cS^{e_1}_1$, we get
 $$\frac{\overline{h}_{p_1}^{\phi_{p_2}}}{\overline{h}_{p_1}} \cdot  (\overline{\eta}^{\phi_{p_2}})^{\frac{p_1-1}{2}}=
 \overline{g}^{\frac{p_1-1}{2}}$$
 in $\overline{\cS^e_k}$. In view of the ``splitting" (\ref{mariaa}) we can derive a contradiction if we check the following:

\begin{lemma}
Assume $\overline{X}$ is the reduction mod $p$ of the open set of the modular curve as in section \ref{modu} and let $\overline{H}=h x^{p-1}$ be the Hasse invariant. Consider the $2$-fold product $\overline{X}^2=\overline{X} \times \overline{X}$ and the projections $\pi_1,\pi_2:\overline{X}^2\ra \overline{X}$. Then there is no rational function $\overline{u}$ such that
\begin{equation}
\label{impos}
\frac{\pi_1^* \overline{h}}{\pi_2^* \overline{h}} =\overline{u}^{\frac{p-1}{2}}.\end{equation}
\end{lemma}

{\it Proof}.
Assume (\ref{impos}) holds for some $\overline{u}$.
Restricting to a horizontal divisor $\overline{X} \times \{\text{point}\}$ we get
$$\overline{h}=\overline{v}^{\frac{p-1}{2}}$$
for some rational function $\overline{v}$ on $\overline{X}$.
Recall that $\overline{H}$ has simple zeroes at the supersingular points; cf. \cite{KM}, 12.4.3. Pick a supersingular point $s$,
let $x_s$ be a basis of the line bundle $L$ in a neighborhood of $s$  and write
$\overline{H}=\overline{h}_s x_s^{p-1}$ with $\overline{h}_s$ having a simple zero at $s$.
Then we get
$$\overline{h}_s=\left( \overline{v}\frac{x^2}{x^2_s}\right)^{\frac{p-1}{2}}.$$
This implies that $\overline{h}_s$ has a zero at $s$ of order divisible by $\frac{p-1}{2}$, a contradiction.
\qed

\bibliographystyle{amsplain}

\end{document}